\documentclass{article}
\usepackage[english]{babel}
\usepackage{amsthm}
\usepackage{amsmath}
\usepackage{amsfonts}
\usepackage{graphicx}
\usepackage{setspace}

\providecommand{\keywords}[1]{\textbf{Key words: } #1}

\def\P{\hbox{\sf P}}

\def\E{\hbox{\sf E}}
\def\M{\hbox{\sf M}}
\def\1{\hbox{\sf 1}}

\title{Stability analysis of two-class  	retrial systems with constant retrial rates and general  service times 
}

\author{Konstantin Avrachenkov  \\	Inria Sophia Antipolis, France\\	\texttt{k.avrachenkov@inria.fr} \\ \and 	Evsey Morozov \\ IAMR KarRC RAS, PetrSU, Petrozavodsk, Russia\\	\texttt{emorozov@karelia.ru}\\
\and 	Ruslana Nekrasova 	 \\ IAMR KarRC RAS, PetrSU, Petrozavodsk, Russia\\	\texttt{ruslana.nekrasova@mail.ru}\\}

\begin{document}
\date{}
\maketitle

\begin{abstract}
We establish stability criterion for a two-class retrial system with Poisson inputs, general class-dependent service times and class-dependent constant retrial rates. We also characterise an interesting phenomenon of partial stability when one orbit is tight but the other orbit goes to infinity in probability. All theoretical results are illustrated by numerical experiments.
\end{abstract}

\begin{spacing}{1.2}
\section{Introduction}

In this work, we consider a two-class retrial system with a single server and no waiting space associated with the server.
If an incoming job finds the server busy, the job goes to the orbit associated with its class. The jobs blocked on a class-dependent orbit attempt to access the server after class-dependent exponential
retrial times in FIFO manner.
The jobs initially arrive to the system
according to Poisson processes. The service times are generally distributed.
The arrival processes as well as service times are class-dependent.

An interested reader can find the description of various types of retrial systems and their applications in the books and surveys \cite{Falin1997,FalinSurvey,Artalejo,Artalejo1,Artalejo2008book}. Specifically, the multi-class retrial systems with constant retrial rate can be applied to computer networks \cite{PEISarticle,Nain}, wireless networks \cite{dim2,dimi,dim3,dimitriou2019power,Dimitriou2020Valuetools,Tuan} and call centers \cite{PhungDuc2016,PEVA}.

Let us outline our contributions and the structure of the paper.
After providing a formal description of the system in Section~\ref{sec:desc}, we first establish equivalence in terms of stability
between the original continuous-time system and a discrete-time system embedded in the departure instants, see Section~\ref{subsec:equiv}. Then, in Section~\ref{subsec:mainres} we prove the stability criterion for our retrial system. In Section~\ref{subsec:balking}
we give an extension of the stability criterion to the modified system with balking, which is useful for modelling two-way communication systems. In Section~\ref{sec:part}
we characterize a very interesting regime of partial stability when one orbit is tight and the other orbit goes to infinity in probability. In particular, we show that in this regime,
as time progresses, the original two-orbit system becomes equivalent to a single orbit system. Curiously enough,
this new single orbit system gains in stability region
due to the jobs lost at infinity. Namely, the stability of one orbit is attained in part  due to a `displacement' of the customers going from other (unstable) orbit,
and it gives a new insight to the transience phenomena.

We mention the most related works in the next paragraph, leaving the detailed
description of related works and the comparison of various stability conditions to Section~\ref{sec:compare}. Finally, in Section~\ref{sec:sim} we demonstrate all theoretical results by simulations with exponential and Pareto distributions of the service times.

Stability conditions for the single-class retrial systems with constant retrial rates have been investigated in \cite{Fayolle1986,Lillo,PEISarticle,KostiaUri,MMOR}.
In \cite{AMNS} the authors have established the necessary stability
conditions for the present system that coincide with the sufficient
conditions obtained here. In fact, the necessary conditions have
been obtained for $N$-orbit systems, with $N \ge 2$. We would like
to note that the proof of the necessary conditions turns out to be
much less challenging than the proof of sufficiency of the same
conditions. In \cite{Nain} the necessary and sufficient conditions
have been established by algebraic methods for the case of two
classes in a completely Markovian setting with the same service
rates. In \cite{dim2} the author, also in the Markovian setting,
has generalized the model of \cite{Nain} to the case of coupled
orbits and different service rates. Then, the author of \cite{Dimitriou2018EJOR}
conjectured sufficient conditions for the two-class retrial
system in the case of general
service times. In \cite{Questa2015}, using an auxiliary majorizing system, the authors have obtained sufficient stability conditions for a very
general multi-class retrial system with $N \ge 2$ classes. Their
sufficient conditions coincide with the necessary conditions of the
present model in the case of homogeneous classes. Recently, the
authors of \cite{PEVA} have also obtained sufficient (but not generally necessary)
conditions for the multi-class retrial system with balking.
To the best of our knowledge, in this paper we for the first time
establish stability criterion for the two-class retrial system
with constant retrial rates and general service times. We credit
this to the combination of the regenerative approach \cite{MSbook2021}
with the Foster-Lyapunov approach for stability analysis of random walks
\cite{FMM}.
Finally, the concept of partial (local) stability has been studied
in \cite{Questa2015} in the context of retrial systems and in \cite{FOSS}
in a more general context of Markov chains. In the present work, we
use both approaches from \cite{Questa2015} and \cite{FOSS} to obtain a refined characterisation of the phenomenon of partial stability in
multi-class retrial systems.

\section{System description}
\label{sec:desc}

Consider a single-server two-class retrial queueing system with constant retrial rates. The system has two Poisson inputs with class-dependent  rates $\lambda_k$ and generic
service times $S^{(k)}$ with distribution functions $F_k$, $k=1,2$. There is no waiting space but two
orbits.
   Define the basic three-dimensional process
\begin{equation}
   \label{cont}
{\bf X}(t)=(N(t),\, X^{(1)}(t),\, X^{(2)}(t)),
\quad t\ge0,
\end{equation}
where $N(t)=1$ if the server is busy at instant $t^-$ ($N(t)=0$ otherwise) and $X^{(k)}(t)$ is the state (size) of orbit $k$ at instant $t^-,\,k=1,2$. If an incoming customer of class-$k$ comes to the system and sees that the server is busy, he/she goes to the $k$-th orbit. The class-$k$ customers retries from orbit-$k$ in FIFO manner with exponential retrial times with rate $\alpha_k$.

In general, the continuous-time process $\{{\bf X}(t),\,t\ge0\}$ is not Markovian.
Now let us construct a discrete-time process, embedded  in the  process $\{{\bf X}(t),\,t\ge0\}$ at the departure instants, which turns out to constitute a Markov chain.
Denote by $\{D_n,\, n\ge 1\}$  the sequence of the departure instants, and   let $X^{(k)}_n=X^{(k)}(D_n)$ be the number of customers in orbit $k$  just after the $n$-th departure, $k=1,\,2$.
Construct the following two-dimensional discrete-time process
\begin{equation}\label{11}
{\bf X}_n = (X^{(1)}_n,\, X^{(2)}_n),\qquad n\ge 1.
\end{equation}
It is easy to check that the
process $\{{\bf X}_n\}$
is a homogeneous irreducible  aperiodic Markov chain (MC).
Let us define  the increments
\begin{equation}
\Delta^{(k)}_{n+1}=X^{(k)}_{n+1}-X^{(k)}_{n},\qquad k=1,\,2,
\end{equation}
and introduce the sequence of vectors
$$
\Delta_{n+1}=(\Delta^{(1)}_{n+1},\, \Delta^{(2)}_{n+1}),\qquad  n\ge 1.
$$
Then, the dynamics of MC $\{{\bf X}_n\}$ dynamics is described by
\begin{equation}
    {\bf  X}_{n+1}={\bf  X}_n+\Delta_{n+1},
\end{equation}
where the distribution of $\Delta_{n+1}$ depends on the value  of ${\bf  X}_n$ only.

\subsection{Transition probabilities of the embedded MC}

Denote by $\mathcal{I}_n$ and $\mathcal{B}_n$, the idle and busy periods of the server between the $n$-th and the $(n+1)$-st departures,  respectively, $n\ge 1$.  Thus, the $(n+1)$-st departure instant can be recursively presented as
$$
D_{n+1}=D_n+\mathcal{I}_n+\mathcal{B}_n.
$$
Then, let $\mathcal{I}$ and $\mathcal{B}$ be the corresponding  generic times.  Next, define by ${\cal A}^{(k)}_n$  the event, when the $(n+1)$-st customer in the server belongs to  class-$k$, $k=1,\,2$.  Note, that on the event ${\cal A}^{(k)}_n$, $\mathcal{B}$ is distributed as service time $S^{(k)}$.
On the other hand, the distribution of $\mathcal{I}$ depends on the state of the orbits: idle/busy.
Now we consider all possible cases separately.

\bigskip

\noindent
1. {\it Both orbits are empty.} In this case ${\bf X}_n=(0,\,0)$
and the server stays idle until the next arrival. Thus, the idle period $\mathcal{I}$
is exponentially distributed with rate $\lambda_1+\lambda_2$ and the mean
$\E \mathcal{I}=1/(\lambda_1+\lambda_2)$.

Denote by $p_k^{(1)}(i)$
the probability that  $i$
customers  join the $1$-st orbit
in the interval $[D_n,\,D_{n+1})$, provided  that both orbits are empty at instant $D_n$ and the $(n+1)$-st customer arriving to the server is from class $k$, and let $p_k^{(2)}(i)$ be the similar probability for the 2-nd orbit. Thus, for $k=1, 2;\,i\ge 0$ we can write
\begin{eqnarray}
p_k^{(1)}(i)&=&\int_0^{\infty}e^{-\lambda_1t}\frac{(\lambda_1t)^i}{i!}dF_k(t),\label{p11}\\
p_k^{(2)}(i)&=&\int_0^{\infty}e^{-\lambda_2t}\frac{(\lambda_2t)^i}{i!}dF_k(t).
\label{p22}
\end{eqnarray}
In fact, the following two subcases are possible:

\begin{enumerate}
    \item[1.1] With probability
\begin{equation}
    \label{toM1-}
\frac{\lambda_1}{\lambda_1+\lambda_2}p_1^{(1)}(i),
\end{equation}
a $1$-st class customer occupies the server and
$i$ customers join the $1$-st orbit, resulting in $X^{(1)}_{n+1}=i$.
Moreover, with probability
\begin{equation}
    \label{toM2-}
\frac{\lambda_1}{\lambda_1+\lambda_2}p_1^{(2)}(j),
\end{equation}
this customer, during the service,
generates $j$ new class-2 orbital customers,
resulting in, $X^{(2)}_{n+1}=j$.
    \item[1.2] With probabilities
\begin{equation}
    \label{toM11-}
\frac{\lambda_2}{\lambda_1+\lambda_2}p_2^{(1)}(i)\quad\text{ and }\quad \frac{\lambda_2}{\lambda_1+\lambda_2}p_2^{(2)}(j),
\end{equation}
a $2$-nd class {\it external} arrival captures  the server and
we obtain $X^{(1)}_{n+1}=i$ and $X^{(2)}_{n+1}=j$, respectively, with the above probabilities.
\end{enumerate}

\bigskip

\noindent
2. {\it Only the $1$-st orbit is empty.} Note that in this case $\E\mathcal{I}=1/(\lambda_1+\lambda_2+\alpha_2)$. Then consider the following cases.
\begin{enumerate}
    \item[2.1] With probability
\begin{equation}
    \label{toM1}
\frac{\lambda_1}{\lambda_1+\lambda_2+\alpha_2}p_1^{(1)}(i),
\end{equation}
a $1$-st class customer occupies the server and
$i$ customers join the $1$-st orbit, resulting in $X^{(1)}_{n+1}=i$.
Moreover, with probability
\begin{equation}
    \label{toM2}
\frac{\lambda_1}{\lambda_1+\lambda_2+\alpha_2}p_1^{(2)}(j),
\end{equation}
this customer, during the service,
generates $j$ new class-2 orbital customers, that is,
$X^{(2)}_{n+1}=X^{(2)}_n+j$.
    \item[2.2] With probabilities
\begin{equation}
    \label{toM11}
\frac{\lambda_2}{\lambda_1+\lambda_2+\alpha_2}p_2^{(1)}(i)
\quad\text{ and }\quad \frac{\lambda_2}{\lambda_1+\lambda_2+\alpha_2}p_2^{(2)}(j),
\end{equation}
a $2$-nd class {\it external} arrival captures  the server and
we obtain $X^{(1)}_{n+1}=i$ and $X^{(2)}_{n+1}=X^{(2)}_n+j$, respectively, $i,j\ge0$.

 \item[2.3]
With probabilities
\begin{equation}
    \label{toM22}
\frac{\alpha_2}{\lambda_1+\lambda_2+\alpha_2}p_2^{(1)}(i)\quad\text{ and }\quad \frac{\alpha_2}{\lambda_1+\lambda_2+\alpha_2}p_2^{(2)}(j),
\end{equation}
an {\it orbital} class-2  customer occupies the server and we obtain
$$
X^{(1)}_{n+1}=i\quad\text{ and }\quad X^{(2)}_{n+1}=X^{(2)}_n-1+j,
$$
respectively, $\,i,j\ge0$.
\end{enumerate}

\bigskip

\noindent
3. {\it Only the  $2$-nd  orbit is empty.}
In this  case $X^{(1)}_n>0,\, X^{(2)}_n=0$, and we have
$\E\mathcal{I}=1/(\lambda_1+\lambda_2+\alpha_1)$. Next we consider  the following three possible cases.
\begin{enumerate}
    \item[3.1] A class-$1$ {\it external} arrival occupies the server, $i$ class-$1$ customers join the 1-st orbit
    with probability
    \begin{equation}\label{13}
        \frac{\lambda_1}{(\lambda_1+\lambda_2+\alpha_1)}p_1^{(1)}(i),
    \end{equation}
    implying $X^{(1)}_{n+1}=X^{(1)}_n+i,\,\,i\ge0$, and, simultaneously,  $j$ class-$2$ customers join the orbit with probability
    \begin{equation}
        \frac{\lambda_1}{(\lambda_1+\lambda_2+\alpha_1)}p_1^{(2)}(j),
    \end{equation}
   implying  $X^{(2)}_{n+1}=j,\,j\ge0.$

    \item[3.2] A retrial attempt from the $1$-st class orbit was successful  (a secondary customer occupies the server before the next {\it external} arrival)
    and $i$ and $j$ customers join class-$1$ and class-$2$ orbits
    with probabilities
     \begin{equation}
        \frac{\alpha_1}{(\lambda_1+\lambda_2+\alpha_1)}p_1^{(1)}(i)\quad
        \mbox{and} \quad \frac{\alpha_1}{(\lambda_1+\lambda_2+\alpha_1)}p_1^{(2)}(j),
    \end{equation}
     respectively. As a result, we obtain
     $$
     X^{(1)}_{n+1}=X^{(1)}_n-1+i, \quad X^{(2)}_{n+1}=j.
     $$

    \item[3.3]  The server becomes busy with the $2$-nd class {\it external} arrival and with probabilities
    \begin{equation}\label{16}
        \frac{\lambda_2}{(\lambda_1+\lambda_2+\alpha_1)}p_1^{(1)}(i)\quad \mbox{and} \quad \frac{\lambda_2}{(\lambda_1+\lambda_2+\alpha_1)}p_1^{(2)}(j),
    \end{equation}
    we have
    $$
    X^{(1)}_{n+1}=X^{(1)}_n+i, \quad X^{(2)}_{n+1}=j, \quad i,j \ge0,
    $$
    respectively.
\end{enumerate}

\bigskip

\noindent
4. {\it Both orbits are busy.}
In this case $X^{(1)}_n>0,\,X^{(2)}_n>0$, and similarly to the above,
we obtain $X^{(1)}_{n+1}=X^{(1)}_n+i,\,i\ge0$, with probability
\begin{equation}\label{17}
    \frac{\lambda_1}{(\lambda_1+\lambda_2+\alpha_1+\alpha_2)}p_1^{(1)}(i)+\frac{\lambda_2+\alpha_2}{(\lambda_1+\lambda_2+\alpha_1+\alpha_2)}p_2^{(1)}(i).
\end{equation}
In the case of a successful class-$1$ retrial attempt, we have
$X^{(1)}_{n+1}=X^{(1)}_n-1+i$
with probability
\begin{equation}
    \frac{\alpha_1}{(\lambda_1+\lambda_2+\alpha_1+\alpha_2)}p_1^{(1)}(i).
\end{equation}
Similarly,
$X^{(2)}_{n+1}=X^{(2)}_n+j, \, j\ge 0$ with probability
\begin{equation}
    \frac{\lambda_2}{(\lambda_1+\lambda_2+\alpha_1+\alpha_2)}p_2^{(2)}(j)+\frac{\lambda_1+\alpha_1}{(\lambda_1+\lambda_2+\alpha_1+\alpha_2)}p_1^{(2)}(j),
\end{equation}
and
$X^{(2)}_{n+1}=X^{(2)}_n-1+j$ with probability
\begin{equation}\label{20}
    \frac{\alpha_2}{(\lambda_1+\lambda_2+\alpha_1+\alpha_2)}p_2^{(2)}(j),
    \quad j \ge 0.
\end{equation}

\section{Stability criterion}
\label{sec:criter}

\subsection{Stability of the embedded MC and underlying continuous-time process}
\label{subsec:equiv}
In this section we establish a connection between the notion of {\it stability (ergodicity)} of the embedded MC introduced in the previous section and the concept of {\it positive recurrence}, which is an analogous notion of {\it stability} for regenerative processes in continuous time.
Although it seems quite intuitive that stability of the embedded MC implies the positive recurrence of the underlying continuous-time process
and vice versa, it is instructive to give a formal proof of this fact.

Recall the definition of the basic three-dimensional process
$$
{\bf X}(t)=(N(t),\, X^{(1)}(t),\,X^{(2)}(t)),\quad t\ge0,
$$
where $N(t)$ is the indicator function of the server occupancy at time instant $t^-$ and $X^{(k)}(t)$ is the size of orbit $k$.
Denote by $t_n$ the arrival instants of the (superposed) Poisson process and let $\hat {\bf X}_n={\bf X}(t_n^-),\,n\ge1.$   We stress that the new hat-notation reflects the fact that the discrete-time process $\{\hat {\bf X}_n\}$  in general is {\it not a Markov chain} and evidently  differs from the original Markov chain $\{{\bf  X}_n\}$  obtained by embedding at the {\it departure instants}.

We recall that the process $\{X(t)\}$ is called {\it regenerative} with {\it regeneration instants} $T_n$ defined recursively as
\begin{eqnarray}
T_{n+1}=\min_i(t_i>T_n: \hat {\bf X}_i={\bf 0}), \quad T_0=0, \quad n\ge 0.\label{84a}
\end{eqnarray}
Note that the equality $ \hat {\bf X}_i={\bf 0}$ is component-wise.
We note that $T_n$  represents the arrival instant of such  a customer which  meets the system  totally idle in the $n$th time.
We assume that the  1st customer arrives in an empty system at instant $t=0$. Such a setting is  called {\it zero initial state} \cite{MorozovDelgado}, and the corresponding regenerative process is called zero-delayed \cite{Asmus}. Denote by $T$  generic {\it regeneration period} (which is distributed as any difference $T_{n+1}-T_n$).   Then the regenerative process is called {\it positive recurrent} if $\E T<\infty$.
Denote  by $\tau$ the generic interarrival (exponential)  time in the superposed Poisson input process, which has  rate $\lambda=\sum\lambda_k$.
Because the input is Poisson, the regeneration period $T$ is {\it non-lattice}. Then,  
the positive recurrence implies the existence of the stationary distribution of the process ${\bf X}(t)$  as $t\to \infty$ and hence the stability of the system \cite{Asmus}. To study stability, it is much more convenient to work with a one-dimensional  process
$$
Z(t)=N(t)+X_1(t)+X_2(t), \quad t\ge0,
$$
counting the total number of customers in the system, which is {\it regenerative with the same regeneration instants} (\ref{84a}).
In the following lemma we establish the equivalence between the stability of the embedded MC and the stability of the original continuous-time process for the case of zero initial state.\\

\noindent
{\bf Lemma 1.}  {\it The zero-initial state Markov chain ${\bf  X}$ is positive recurrent  if and only if the process $\{X(t)\}$  is positive recurrent, that is if $\E T<\infty$.}

\begin{proof}
If the process $\{X(t)\}$ is positive recurrent,
then it follows by a regenerative argument \cite{Asmus} that the {\it stationary probability $\P_0$}, the probability of the system being totally free, exists  and is equal to
\begin{eqnarray}
\P_0=\lim_{t\to \infty}\P(Z(t)=0)
=\lim_{t\to\infty}\frac{1}{t}
\int_0^t \1(Z(u)=0)du=\frac{\E \tau}{\E T}>0, \quad \mbox{w.p. 1,} \label{84}
\end{eqnarray}
where $\tau$  denotes the generic inter-arrival time in the superposed Poisson input process.
On the other hand,
$\{\hat{\bf X}_n\}$ is the  embedded  discrete-time   regenerative process with the regeneration instants
\begin{eqnarray}
\theta_{n+1}=\min(i: i>\theta_n,\, \hat {\bf X}_i={\bf 0}),\quad \theta_0=0, \quad n\ge 0,\label{77a}
\end{eqnarray}
and $\theta$  denoting generic regeneration period of this  discrete-time process. Namely, the generic regeneration cycle is given by
$\theta =_{st} \theta_{n+1}-\theta_n$.
It is well known that the discrete-time length $\theta$ of the regeneration cycle
is connected with the continuous-time
length $T$ by the following stochastic equality \cite{Asmus}
(Chapter X, Propositions 3.1 and 3.2):
$$
T=_{st}\sum_{k=1}^\theta\tau_k\,\,\,(\sum_\emptyset:=0),
$$
where $\tau_k=t_{k+1}-t_k$ is the $k$-th inter-arrival interval and the summation index $\theta$ is a (randomized) stopping time.
It then immediately follows from the Wald's identity that
$$
\E T=\E\tau \E \theta.
$$
Note that, because $0<\E\tau<\infty$, then $\E T=\infty$ implies $\E \theta=\infty$, and vice versa.
Thus we obtain that  $\E \theta<\infty$, and hence the positive recurrence of the basic process $\{{\bf X}(t)\}$ implies the positive recurrence of
the process  $\{\hat {\bf X}_n\}$  embedded (in the basic process)
at the arrival instants.
Conversely, $\E \theta<\infty$ implies $\E T<\infty$ as well.

It remains to connect the process $\{\hat {\bf X}_n\}$ with the embedded MC we studied above.
Note that, as in (\ref{84a}), regenerations $\{\theta_n\}$,  defined in (\ref{77a}), are generated by the arrivals meeting empty system.
On the other hand,  $\theta$ represents {\it both the number of arrivals and the number of departures} from the system within a {\it continuous-time regeneration period} $T$.
Thus, $\theta$  is also generic regeneration period of the embedded MC $\{{\bf  X}_n\}$.
It is worth mentioning that, at each instant of time, the index of a customer which see empty system (and generating new regeneration period of the  processes $\{{\bf X}(t)\}$  and $\{\hat {\bf X}_n\}$),  differs from the index of a customer {\it leaving empty system} by no more than one.
Thus we obtain that
$$
\hat \pi_0=\lim_{n\to \infty}\P(\hat {\bf X}_n={\bf 0})=\lim_{n\to\infty}\P({\bf X}_n={\bf 0})
=:\pi_0=\frac{1}{\E \theta}>0.
$$
 Hence
$\pi_0>0$, and because
the MC we consider  is aperiodic and irreducible,
then it is also ergodic.
Thus we see that  both concepts of stability (in continuous and discrete time) agree and the lemma hereby is proven.
\end{proof}

We note that, at the first sight, the equality $\pi_0={\hat \pi}_0$ seems rather surprising because $\pi_0$ relates to the MC while ${\hat \pi}_0$  relates to the regenerative process $\{\hat{\bf  X}_n\}$ which in general is not Markov.

\subsection{Stability  of the embedded Markov chain}
\label{subsec:mainres}

The results of this section are based on the general stability conditions for two-dimensional MCs, obtained in \cite{FMM}.

Let us first introduce some additional notations for the embedded Markov chain (\ref{11}) needed for the
application of the results from \cite{FMM}. Specifically, let us derive the drifts of the embedded MC in various regions of the state space.

Denote by $\M^{01}_k$ the  mean number of customers joining the class-$k$ orbit in the time interval $[D_n,\,D_{n+1})$, provided  $X^{(1)}_n=0,\, X^{(2)}_n>0$.
Recall that $\mu_k=1/\E S^{(k)}$ and denote
$$
\rho_k=\frac{\lambda_k}{\mu_k},\,\,\,\,\,\hat \rho_k=\frac{\alpha_k}{\mu_k},\,\,\,k=1,\,2.
$$
Then it follows from  (\ref{toM1})--(\ref{toM22}) after a simple algebra that
\begin{eqnarray}\label{m01}
\M^{01}_1&=&\sum_{i\ge 1}i\Bigl[\frac{\lambda_1p_1^{(1)}(i)}{\lambda_1+\lambda_2+\alpha_2}+\frac{(\lambda_2+\alpha_2)p_2^{(1)}(i)}{\lambda_1+\lambda_2+\alpha_2}\Bigr]=
\frac{\lambda_1(\rho_1+\rho_2+\hat\rho_2)}
{\lambda_1+\lambda_2+\alpha_2},\\
\M^{01}_2&=&\sum_{i\ge 1}i\Bigl[\frac{\lambda_1p_1^{(2)}(i)}{\lambda_1+\lambda_2+\alpha_2}+\frac{\lambda_2p_2^{(2)}(i)}{\lambda_1+\lambda_2+\alpha_2}+\Bigr]+\sum_{i\ge 2}\frac{(i-1)\alpha_2p_2^{(2)}(i)}{\lambda_1+\lambda_2+\alpha_2}\nonumber\\
&=&\sum_{i\ge 1}i\Bigl[\frac{\lambda_1p_1^{(2)}(i)+\lambda_2p_2^{(2)}(i)+\alpha_2p_2^{(2)}(i+1)}{\lambda_1+\lambda_2+\alpha_2}\Bigr]\nonumber\\
&=&\frac{\lambda_2(\rho_1+\rho_2+\hat\rho_2)
-\alpha_2} {\lambda_1+\lambda_2+\alpha_2}.\label{m02}
\end{eqnarray}
Similarly, denote by $\M^{10}_k$ the  mean number of customers joining the class-$k$ orbit in the time interval $[D_n,\,D_{n+1})$,  provided  $X^{(1)}_n>0,\, X^{(2)}_n=0,\,k=1,2$. 
Then by analogy with (\ref{m01}) and (\ref{m02})  from (\ref{13})--(\ref{16}) we  obtain 
\begin{eqnarray}\label{beMM1}
\M^{10}_1&=&\frac{\lambda_1(\rho_1+\rho_2+\hat \rho_1)
-\alpha_1} {\lambda_1+\lambda_2+\alpha_1},\\
\M^{10}_2&=&\frac{\lambda_2(\rho_1+\rho_2+\hat\rho_1)} {\lambda_1+\lambda_2+\alpha_1}.\label{beMM}
\end{eqnarray}
Continuing in the same way, we denote by  $\M^{11}_k$ the mean number of customers joining the class-$k$ orbit in the time interval $[D_n,\,D_{n+1})$, if  $X^{(1)}_n>0,\, X^{(2)}_n>0$, and  from (\ref{17})--(\ref{20}) we obtain 
\begin{eqnarray}\label{beM}
\M_1^{11}
&=&\frac{(\lambda_1+\alpha_1)\rho_1+
\lambda_1(\rho_2+\hat\rho_2)-\alpha_1}{\lambda_1+\lambda_2+\alpha_1+\alpha_2},
\end{eqnarray}
\begin{eqnarray}\label{beM2}
\M_2^{11}
&=&\frac{(\lambda_2+\alpha_2)\rho_2+\lambda_2(\rho_1+\hat\rho_1)-\alpha_2}{\lambda_1+\lambda_2+\alpha_1+\alpha_2}.
\end{eqnarray}

Our further analysis is based on Theorem A presented in the Appendix, a statement from \cite{FMM}.
Note, that in the general case,  Theorem A is applicable  under some additional technical conditions (see Appendix A),
which hold automatically when the input is Poisson.
Denote the total load coefficients by
\begin{eqnarray*}
\rho=\rho_1+\rho_2,\,\,\,
\hat\rho=\hat \rho_1+\hat \rho_2.
\end{eqnarray*}
Now we are in a position to state our central result.

\bigskip

\noindent
{\bf Theorem 1.}
{\it Two-class retrial system with constant retrial rates, Poisson inputs, general service times and exponential retrials is
ergodic if and only if
\begin{equation}\label{Theo2}
 \rho<\min \Big( \frac{\hat\rho_1}{\rho_1+\hat\rho_1},\,
 \frac{\hat\rho_2}{\rho_2+\hat\rho_2}
 \Big).
\end{equation}
}

\begin{proof}
Note that the defined above drifts $\M^{01}_{k},\,\M^{10}_{k},\,\M^{11}_{k},\,k=1,2$,
correspond to the drifts used in the statement of Theorem~A (see the Appendix).

First, let us consider the conditions for the case (a) of Theorem A.
Specifically, the condition $\M^{11}_1<0$
takes the following form
\begin{eqnarray}\label{IdM1}
\rho_1\Big(\lambda_1+\alpha_1+(\lambda_2+\alpha_2)\frac{\mu_1}{\mu_2}\Big)&<&\alpha_1,
\end{eqnarray}
while the condition $\M^{11}_2<0$ takes the form
\begin{eqnarray}
\rho_2\Big(\lambda_2+\alpha_2+(\lambda_1+\alpha_1)\frac{\mu_2}{\mu_1}\Big)&<&\alpha_2.\label{IdM2}
\end{eqnarray}
The inequalities (\ref{IdM1}) and (\ref{IdM2}) can be further rewritten as
\begin{eqnarray*}
\rho_k\Big(\rho+\hat\rho\Big)<\hat\rho_k,\,\,k=1,\,2.
\end{eqnarray*}
Next, the first condition $\M_1^{11}\M_2^{10}-\M_2^{11}\M_1^{10}<0$
in (\ref{first}) becomes, after a tedious algebra (see Appendix B for details),
\begin{equation}\label{alll1}
\rho<\frac{\hat\rho_1}{\rho_1+\hat\rho_1},
\end{equation}
while the second condition $\M_2^{11}\M_1^{01}-\M_1^{11}\M_2^{01}<0$ in (\ref{first})
can be transformed to
\begin{equation}\label{alll2}
\rho<\frac{\hat\rho_2}{\rho_2+\hat\rho_2}.
\end{equation}

Now we can express in terms of load coefficients, the three ergodic cases (a.1), (b.1) and (c.1) of Theorem~A (see inequalities (\ref{first}),  (\ref{first_b}) and (\ref{first_c}) in Appendix~A) as follows:

\bigskip

\noindent
Case (a.1)
    \begin{equation}\label{aa}
\left\{
\begin{gathered}
\rho_k\Big(\rho+\hat\rho\Big)<\hat\rho_k, \,\,\,k=1,2,\hfill \\
\rho<\frac{\hat\rho_k}{\rho_k+\hat\rho_k}
,
\quad k=1,\,2; \hfill
\\
\end{gathered}
\right.
\end{equation}
 Case (b.1)
    \begin{equation}\label{bb}
\left\{
\begin{gathered}
\rho_1\Big(\rho+\hat\rho\Big)\ge\hat\rho_1, \hfill \\
\rho_2\Big(\rho+\hat\rho\Big)<\hat\rho_2,\hfill \\
\rho<\frac{\hat\rho_1}{\rho_1+\hat\rho_1};
\hfill \\
\end{gathered}
\right.
\end{equation}
 Case (c.1)
    \begin{equation}\label{cc}
\left\{
\begin{gathered}
\rho_1\Big(\rho+\hat\rho\Big)<\hat\rho_1, \hfill \\
\rho_2\Big(\rho+\hat\rho\Big)\ge \hat\rho_2,\hfill \\
\rho<\frac{\hat\rho_2}{\rho_2+\hat\rho_2}.
\hfill
\\
\end{gathered}
\right.
\end{equation}

Now our goal is to simplify ergodicity
conditions (\ref{aa})--(\ref{cc}).
Towards this goal, we rewrite the system (\ref{aa})
in terms of functions of $\alpha_1$ and $\alpha_2$, assuming other parameters fixed.
The first pair of inequalities in (\ref{aa}) can be transformed to
\begin{eqnarray}\label{g1}
\alpha_2&<&\frac{1-\rho_1}{\rho_1}\frac{\mu_2}{\mu_1}\, \alpha_1-\rho \mu_2=:g_1(\alpha_1),\\
\alpha_2&>&\frac{\rho_2}{1-\rho_2}\frac{\mu_2}{\mu_1} \, \alpha_1+\frac{\rho_2}{1-\rho_2}\rho \mu_2=:g_2(\alpha_1),
\label{g2}\end{eqnarray}
while the pair of inequalities
$\rho<\alpha_k/(\lambda_k+\alpha_k), k=1,2,$ becomes
\begin{eqnarray}
\alpha_k&>&\frac{\rho}{1-\rho}\lambda_k,\quad k=1,2.
\label{G2}
\end{eqnarray}
For fixed values $\lambda_k,\,\mu_k$, such that $\rho<1$,
the right hand sides of (\ref{g1}), (\ref{g2})
are the increasing linear functions of $\alpha_1$ with
a common point $(\alpha_1^*,\,\alpha_2^*)$, where
\begin{eqnarray}\label{G-1}
 \alpha_k^*=\frac{\rho}{1-\rho}\lambda_k,\quad k=1,2.
  \end{eqnarray}
The ergodic case (a.1), described by system (\ref{aa}),
corresponds to the values of $(\alpha_1,\,\alpha_2)$ such that
\begin{eqnarray}\label{forPr}
& g_2(\alpha_1) < \alpha_2 < g_1(\alpha_1),& \\
& \alpha_1 > \alpha_1^*. &
\label{forPr2}\end{eqnarray}
Let us now show that the set of $(\alpha_1,\,\alpha_2)$ corresponding to (\ref{forPr}) and (\ref{forPr2}) is non-empty as long as $\rho < 1$. Assume on contrary that
\begin{equation}\label{asLin}
  g_2(\alpha_1)\ge g_1(\alpha_1),
\end{equation}
thus (\ref{forPr}) is violated. The inequality (\ref{asLin}) for linear increasing functions $g_1,\,g_2$ under conditions $\alpha_k>\alpha_k^*,\,k=1,\,2$ implies a similar relation for the coefficients in front of $\alpha_1$ (see (\ref{g1}),(\ref{g2})), that is
\begin{equation}\label{coll1}
    \frac{\rho_2}{1-\rho_2}\frac{\mu_2}{\mu_1}\ge \frac{1-\rho_1}{\rho_1}\frac{\mu_2}{\mu_1}.
\end{equation}
Multiplying both parts of (\ref{coll1}) by $\rho_1(1-\rho_2)\mu_1/\mu_2\ge 0$, we obtain
$$
\rho_1\rho_2\ge (1-\rho_1)(1-\rho_2),
$$
which is equivalent to $\rho = \rho_1+\rho_2 > 1$ and yields a contradiction.

Now, similarly, we describe the stability regions (b.1) and (c.1), presented in (\ref{bb}) and (\ref{cc}), in terms of the functions
$g_1(\alpha_1)$ and $g_2(\alpha_1)$ as follows:
\begin{eqnarray}
    \mbox{Case (b.1): } & & \alpha_2\ge
    g_1(\alpha_1),\quad \alpha_1>\alpha_1^*,\nonumber \\
    \mbox{Case (c.1): } & & \alpha_2\le
    g_2(\alpha_1),\quad \alpha_2>\alpha^*_2.\nonumber
\end{eqnarray}

Next, by combining the three cases, we conclude that the embedded MC is ergodic, if $\alpha_1>\alpha_1^*$ and $\alpha_2>\alpha_2^*$, which is equivalent in fact to $\rho<\alpha_k/(\lambda_k+\alpha_k),\,k=1,2$.
Thus, the conditions (\ref{aa})--(\ref{cc}) can be written as (\ref{Theo2}).

Moreover, we can also delimit  the transience regions  in terms of $\alpha_k,\,g_k(\alpha_1),\,k=1,\,2$ (see Theorem~A in the Appendix).
Figure \ref{fig1} illustrates stability/transience regions in different cases for a fixed $\rho<1$.

\begin{figure}[htb]	
\centering
\includegraphics[width=0.8\textwidth]{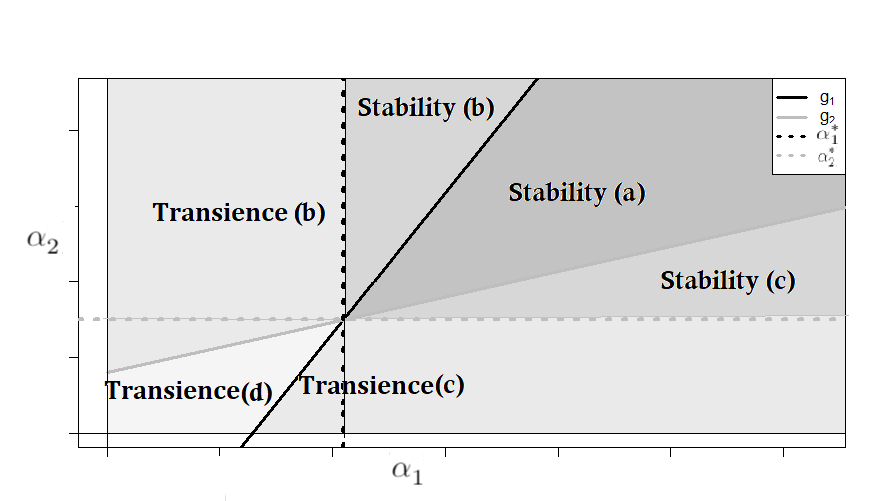}
	\caption{Stability/transience regions.
	}
	\label{fig1}
\end{figure}

Note,  if  (\ref{Theo2}) is violated, then the basic MC $\{{\bf  X}_n\}$ is transient
by Theorem A.   In this case, using the proof by contradiction and  regenerative approach, one can show that  at least one component of this two-dimensional vector goes to infinity in probability, see for instance,  \cite{MorozovDelgado}.

Thus
(\ref{Theo2}) is a sufficient stability (ergodicity) condition.
To show that this condition is also necessary, we refer to the paper \cite{AMNS}  where it is shown, in the  adopted notation and with $\rho=\rho_1+\cdots+\rho_N$,
that if $N$-class retrial system with Poisson  inputs is
ergodic
then
\begin{equation}\label{forN}
    \lambda_k\rho<\alpha_k\big(1-\rho\big),\quad k=1,\dots,N.
\end{equation}
(Indeed, in  paper \cite{AMNS},  we apply  an equivalent notion positive recurrence
in the framework of the  regenerative approach, see Lemma 1 above.)
We can rewrite (\ref{forN}) (for $N$=2) 
as
 $$
 \lambda_k\rho\equiv \lambda_k\rho<\alpha_k\big(1-\rho\big),\quad k=1,\,2,
 $$
 implying
\begin{equation}\label{for2}
    \rho<\min_{k=1,2} \frac{\alpha_k}{\lambda_k+\alpha_k}.
\end{equation}
Thus (\ref{for2}) coincides with (\ref{Theo2}) and 
condition (\ref{Theo2})  is also the necessary stability condition.
\end{proof}

{\bf Remark.} {\it It follows from (\ref{for2}) that
if the two-class retrial system under consideration is  ergodic then
$
\rho<1.
$
}

\bigskip

\subsection{Stability of a system with balking}
\label{subsec:balking}

We can assume an extra feature in the system under consideration as follows. If a primary class-$k$ customer meets busy server, he  joins the corresponding orbit with a given ({\it balking}) probability $b_k\ge 0,\,k=1,\,2$ and leaves the system with probability $1-b_k$.
In this case, the stability condition of Theorem~1 transforms to
\begin{equation}\label{balk}
    b_1\rho_1+b_2\rho_2<\min \Big(\frac{\alpha_1}{\alpha_1+\lambda_1},\, \frac{\alpha_2}{\alpha_2+\lambda_2}\Big).
\end{equation}
This is an immediate extension. Namely, taking into account
balking policy, we redefine the transition probabilities (\ref{p11}), (\ref{p22}), and the statement (\ref{balk}) is
then proved by the same arguments as Theorem~1. We note that this modification of the system is useful to model two-way communication systems, for more details see e.g., \cite{Tuan,PEVA}.

\section{Partial stability}
\label{sec:part}

Let us now discuss an effect of {\it partial  stability} to the best of our knowledge first discovered in \cite{Questa2015}.
In the case of two classes of customers, the statement of Theorem~4 from \cite{Questa2015} can be qualitatively formulated as follows: under some (given  below) conditions, class-1 orbit size stays tight while class-2 orbit increases unlimitedly in probability.
(Of course, by the symmetry, this can be formulated for the opposite case, when the 2nd orbit is tight while the 1st orbit increases.)
By the evident reason, this statement can be regarded as the case of {\it partial stability}.

The purpose of this section is firstly to show that, in terms of the embedded MC $\{{\bf X}_n\}$, the partial stability corresponds to
the transient case (c.2) of Theorem~A, i.e., $\M_1^{11}<0,\,\M_2^{11} \ge 0$ and condition $\M_2^{11}\M_1^{01}-\M_1^{11}\M_2^{01}>0$.

Secondly, by establishing a relation with a single-orbit system, we shall show how to describe the long run behaviour of the stable orbit.

Note that the stability conditions which correspond to transience case (c.2) can be defined in terms of the load coefficients as follows:
\begin{eqnarray}
 \hat\rho_1&>&\rho_1(\rho+\hat\rho),  \label{part1} \\
 \hat\rho_2  &\le&   \rho_2(\rho+\hat\rho),\label{part2}\\
\label{part3}
\rho&>&\frac{\hat\rho_2}{\rho_2+\hat\rho_2} .
   \end{eqnarray}
 It is important to  note that  (\ref{part1})  can be  written as
\begin{eqnarray}\label{64}
\rho_1<\frac{\hat\rho_1}{\rho+\hat\rho_1+\hat\rho_2}<1.
\end{eqnarray}
Now we  show that, provided  conditions (\ref{part1}) and  (\ref{part3})   hold, then they imply
condition  (\ref{part2}), which  turns out to be {\it redundant}.

Now we consider  in detail  inequalities  (\ref{part1})--(\ref{part3})  in all three possible sub-cases when  $\rho_1<1$.\\

\noindent
{\it  Sub-case 1: $\rho_1<1,\,\rho<1$.} 
In this case it is convenient to rewrite conditions (\ref{part1})--(\ref{part3})  as  follows:
\begin{eqnarray}\label{case11-1}
\alpha_1&>&\frac{\rho_1}{(1-\rho_1)}\frac{\mu_1}{\mu_2} \alpha_2+\frac{\rho_1}{1-\rho_1}
\lambda_1,
\\
\alpha_1&\ge&\frac{(1-\rho_2)}{\rho_2}\frac{\mu_1}{\mu_2}\alpha_2-\lambda_1
, \label{case11-2}
\\\label{case11-3}
\alpha_2&<& \frac{\rho}{1-\rho}\lambda_2=: \alpha_2^*,
\end{eqnarray}
respectively. Next  assume that the following relation holds  between the  r.h.s. of conditions (\ref{case11-1}) and (\ref{case11-2})
\begin{equation}\label{rh}
    \frac{\rho_1}{(1-\rho_1)}\frac{\mu_1}{\mu_2} \alpha_2+\frac{\rho_1}{1-\rho_1}\lambda_1\le \frac{(1-\rho_2)}{\rho_2}\frac{\mu_1}{\mu_2}\alpha_2-\lambda_1.
\end{equation}
After some algebra, (\ref{rh}) transforms to the inequality
$$
\alpha_2\ge \frac{\rho}{1-\rho}\lambda_2,
$$
which  contradicts (\ref{case11-3}). Thus, we have
\begin{eqnarray}
 \frac{\rho_1}{(1-\rho_1)}\frac{\mu_1}{\mu_2} \alpha_2+\frac{\rho_1}{1-\rho_1}\lambda_1> \frac{(1-\rho_2)}{\rho_2}\frac{\mu_1}{\mu_2}\alpha_2-\lambda_1,\label{69}
\end{eqnarray}
and  inequality  (\ref{case11-1})  implies inequality (\ref{case11-2}); or equivalently, (\ref{part1}) implies (\ref{part2}).  Thus, the latter condition is redundant.

\medskip

\noindent{\it Sub-case 2: $\rho_1<1,\,\rho_2>1$.}
In  this case, by rewriting condition (\ref{part2})
as
\begin{eqnarray}
\label{case2-2}
\rho_2(\rho+\hat\rho_1)\ge \hat\rho_2(1-\rho_2),
\end{eqnarray}
we see that, since the right hand side is negative,
the condition (\ref{part2}) always holds and
hence is redundant in this sub-case.

\medskip

\noindent
{\it Sub-case 3: $\rho_1<1,\,\rho_2<1,\,\rho>1$.}
In this case  conditions (\ref{case11-1}), (\ref{case11-2}) remain unchanged,
while condition (\ref{case11-3})
becomes
\begin{eqnarray}\label{case31-3}
\alpha_2&>& \frac{\rho}{1-\rho}\lambda_2.
\end{eqnarray}
As in sub-case 1, it is easy to check that, provided inequality (\ref{case31-3}) holds,
then condition (\ref{69}) holds as well, and thus condition (\ref{case11-2}); or equivalently, condition  (\ref{part2}) is redundant again.

Consequently a pair of conditions $\M_1^{11}< 0$ and  $\M_2^{11}\M_1^{01}-\M_1^{11}\M_2^{01}>0$ (or its analogues (\ref{part1}) and (\ref{part3})) define the transience case (c.2).
Before we formulate our next results, let us
recall the definition of the failure rate
$$
r(x):=\frac{f(x)}{1-F(x)},
$$
of a non-negative absolutely continuous distribution $F$ with density $f$, defined for all $x$ such that
$1-F(x)>0$. We say that a distribution belongs to class $\mathcal{D}$ if its failure rate
satisfies $\inf_{x\ge 0} r(x)>0$. (Some, fairly common,   distributions satisfying this requirement
can be found  in  \cite{Questa2015}.)

\bigskip

\noindent{\bf Theorem 2.}
{\it If, in  the initially empty system,  conditions
\begin{eqnarray}
\label{61}
 \hat\rho_1>\rho_1(\rho+\hat\rho),\quad
\rho>\hat\rho_2/(\rho_2+\hat\rho_2),
\end{eqnarray}
hold and distribution $F_k$ of
service time of class-k customers  belongs to
class $\mathcal{D},\,k=1,2,$
then the 1-st orbit is tight and the 2-nd orbit increases in probability,
that is $X^{(2)}(t)\Rightarrow \infty$.}

\begin{proof}
Recall notation
$$
\rho=\rho_1+\rho_2,\,\,\,\hat\rho_k=\frac{\alpha_k}{\mu_k},\,\,\,\hat\rho=\hat\rho_1+\hat\rho_2,
$$
and also  denote
\begin{eqnarray}
    \P_L=\frac{\rho+\hat\rho}{1+\rho+\hat\rho}.\label{55}
\end{eqnarray}
Following \cite{Questa2015}, we consider  an auxiliary two-class  system with  two Poisson inputs  with rates $\tilde\lambda_k=\lambda_k+\alpha_k,\,k=1,2$ and the same service times as in the  original system.
In this new system, any class-$k$ customer meeting server busy
becomes ``colored'' and  joins a  virtual orbit being a part of an {\it infinite} class-$k$ queue, which in turn is a source of the Poisson input with rate $\alpha_k$. (For details see  \cite{Questa2015}.)  Then $P_L$ is the stationary ``loss'' probability in this auxiliary system,  that is  the probability that a customer meets server busy.
It  easily  follows from  \cite{Questa2015} that   class-k orbit size (the number of colored customers) in the  auxiliary  system stochastically dominates class-k orbit size in  the  original system, provided  $F_k\in \mathcal{D},\,k=1,2$.
Moreover, it is  shown in \cite{Questa2015} (Theorem 4 there)  
that, if the system is initially empty, and the following conditions  hold:
\begin{eqnarray}
    \P_L&<&\frac{\hat\rho_1}{\hat\rho_1 +\rho_1},\label{57}\\
    \rho &\ge &\frac{\hat\rho_2}{\hat\rho_2 +\rho_2},\label{56}
\end{eqnarray}
then the 1st orbit is tight and $X^{(2)}(t)\Rightarrow \infty$.
\\

\noindent{\bf Remark.} The proof of tightness in \cite{Questa2015}  is based on a monotonicity property which in turn in general has been  proved  for  $F_k\in \mathcal{D}$ only. We believe that this requirement is only technical and indeed is not needed for stability. This conjecture, in particular, is  supported in Section 6.2 by a numerical example with Pareto service time which does not belong to class  $\mathcal{D}$.\\

Thus, our goal is to  show that  conditions (\ref{61}) of Theorem~2 coincide with conditions (\ref{57}), (\ref{56}), and, for this purpose,
we write
conditions (\ref{61}) separately as
\begin{eqnarray}
\hat\rho_1&>&\rho_1(\rho+\hat\rho),  \label{59}\\  \label{591}
\rho&>&\frac{\hat\rho_2}{\rho_2+\hat\rho_2},
   \end{eqnarray}
   respectively.
   Because
$$
\P_L=\frac{\rho+\hat\rho}{1+\hat\rho+\rho}<\frac{\hat\rho_1}{\hat\rho_1+\rho_1},
$$
or, equivalently,
\begin{equation}\label{for_t3}
\rho\rho_1
+\hat\rho\rho_1
<\hat\rho_1,
\end{equation}
we see that (\ref{57}) coincides with (\ref{59}).
It remains to note  that
(\ref{591}) is a particular case of condition (\ref{56}).
Thus, conditions (\ref{61})  define transience case (c.2)
and simultaneously are the assumptions of Theorem 2.
Hence, transience case (c.2) means that the 1st orbit is tight and $X^{(2)}(t)\Rightarrow \infty$.
\end{proof}

\bigskip

Theorem 2 (as well as  Theorem 4 in  \cite{Questa2015})  shows that, unlike in classical retrial systems,  in the constant retrial rate system, stability/instability may happen {\it locally}.
Denote
  \begin{eqnarray}
   \P_B^{(2)}=
   \frac{\rho+\hat\rho_2}{\rho_2+\hat\rho_2+1}.\label{54}
  \end{eqnarray}
 It is shown in \cite{Questa2015}, under conditions (\ref{56}), (\ref{57}),  that $\P_B^{(2)}$ is the limiting fraction of the mean busy time of server, that is, in an evident notation,
$$
\lim_{t\to\infty}\frac{\E B(t)}{t}=\P_B^{(2)}.
$$
Next, we note that the condition (\ref{56}) can be written as
\begin{eqnarray}
\frac{\rho+\hat\rho_2}{\rho_2+\hat\rho_2+1}=\P_B^{(2)}\le \rho.
 \label{77}
 \end{eqnarray}
In particular, if inequality (\ref{56})  is strict, then $\P_B^{(2)}< \rho$. On the other hand,  the equality in (\ref{77}) means that the limiting fraction of busy time $\P_B^{(2)}$ is equal to $\rho$.
This surprising   result  has the following intuitive explanation
(first remarked in \cite{AMNS}):  under strict  inequality (\ref{56}),
a {\it non-negligible} fraction of class-2 customers, joining an infinitely increasing orbit 2, in fact `disappears' from the system, and thus the limiting  fraction of the `processed' work {\it becomes less than}  $\rho$, an arrived  work per unit of time.  If
the equality in (\ref{56}) holds, then the 2nd orbit size increases unlimitedly in general {\it in probability} only.
In this case, the fraction of class-2 customers, joining an infinitely increasing orbit
turns out to be {\it negligible}. As a result, the limiting fraction of the mean busy time (`busy probability'  $\P_B^{(2)}$)  coincides with
the stationary busy probability in the  positive recurrent  case.

Again assume that the conditions (\ref{56}) and  (\ref{57})  hold, and that
$\P_B^{(2)}\ge\P_L$, that is
 $$
 (\rho+\hat\rho_2)(1+\hat\rho+\rho)\ge (\rho+\hat\rho)(\rho_2+\hat\rho_2+1).
 $$
 After a simple algebra, we obtain the inequality
 $$
 \rho_1(\rho+\hat\rho)\ge \hat\rho_1,
 $$
 contradicting (\ref{57}).  Thus, $\P_B^{(2)}<\P_L$.
This result  has
the following intuitive explanation  in the bufferless setting.
Note that the loss  probability $\P_L$ in the auxiliary system,
by PASTA, is also the limiting fraction of server busy time.  Because the 1st orbit in the original system is {\it tight} then it follows that the  fraction of the server idle time  in the original system   is  non-negligible, and as  a result,  this (limiting)  fraction $\P_B^{(2)}$  is strictly  dominated by  the  probability $\P_L$.

Note that the simulation results in Section~6 illustrate the phenomenon of partial stability in regions 6 and 8 (see Figure~\ref{sim1} below).

\subsection{Relation to single-orbit system}

In this section, we establish an intuitively  expected result that the stability conditions from Theorem 2 coincide with stability conditions of the following associated single-server two-class  system: while class-1 customers meeting server busy join the orbit, class-2 customers arrive  as if the 2nd orbit would be permanently busy.
In other words, the 2nd orbit is the  source  of the Poisson process with rate $\alpha_2$.  Evidently, the associated system can be considered as a `limit' of the original system under conditions of Theorem 2. Thus, under those conditions, class-2 customers arrive
to the server with the input rate
$$
\tilde\lambda_2=\lambda_2+\alpha_2,
$$
and leave the system if they find the server busy. While
the external class-1 customers arrive to the server with
the input rate $\lambda_1$.
In this single-orbit system, we denote by $\{X(t)\}$  the orbit-size process, and  let $X_n$ be  the orbit size just after the $n$th departure, $n\ge1$. Recall notation
$p_1^{(1)}(i),\,p_2^{(1)}(i)$ from
(\ref{p11}), (\ref{p22}) and note
that, if
$X_n=0$, then
$X_{n+1}=i$
with the  probability
\begin{equation}
    \frac{\lambda_1}{\lambda_1+\tilde \lambda_2}p_1^{(1)}(i)+\frac{\tilde \lambda_2}{\lambda_1+\tilde \lambda_2}p_2^{(1)}(i),\,\,i\ge0.
\end{equation}
If $X_n>0$, then $X_{n+1}=X_n+i$ with probability
$$
\frac{\lambda_1}{\lambda_1+\alpha_1+\tilde \lambda_2}p_1^{(1)}(i)+\frac{\tilde \lambda_2}{\lambda_1+\alpha_1+\tilde \lambda_2}p_2^{(1)}(i),
$$
and  $X_{n+1}=X_n+i-1$ with probability
$$
\frac{\alpha_1}{\lambda_1+\alpha_1+\tilde \lambda_2}.
$$
This gives  the (conditional) mean orbit size:
\begin{eqnarray*}
&&\E\big[X_{n+1}|X_{n}=i\big]=\sum_{j\ge 0}\frac{\alpha_1}{\lambda_1+\alpha_1+\tilde \lambda_2}p_1^{(1)}(i)\big(i-1+j\big)\\
&+&\sum_{j\ge 0}\Big( \frac{\lambda_1}{\lambda_1+\alpha_1+\tilde \lambda_2}p_1^{(1)}(i)+ \frac{\tilde\lambda_2}{\lambda_1+\alpha_1+\tilde \lambda_2}p_2^{(1)}(i)
\Big)\\
&=&i+\frac{\lambda_1+\alpha_1}{\lambda_1+\alpha_1+\tilde \lambda_2}\rho_1+\frac{\tilde\lambda_2}{\lambda_1+\alpha_1+\tilde \lambda_2}\frac{\lambda_1}{\lambda_2}\rho_2-\frac{\alpha_1}{\lambda_1+\alpha_1+\tilde \lambda_2}.
\end{eqnarray*}
Thus,  the negative drift
$
\E\big[X_{n+1}|X_{n}=i\big]<i
$
holds if
$$
\frac{\lambda_1+\alpha_1}{\lambda_1+\alpha_1+\tilde \lambda_2}\rho_1+\frac{\tilde\lambda_2}{\lambda_1+\alpha_1+\tilde \lambda_2}\frac{\lambda_1}{\lambda_2}\rho_2<\frac{\alpha_1}{\lambda_1+\alpha_1+\tilde \lambda_2}.
$$
The latter inequality, after a simple algebra, becomes
\begin{equation}
    \label{next}
\rho_1(\rho+\rho_1)<\hat\rho_1,
\end{equation}
and  coincides with condition (\ref{59})  implying tightness of the 1st orbit in partially stable scenario for the two-orbit system. 

Next we  establish a stronger result: the weak {\it convergence} of the two-orbit system to the associated  single-orbit system.
To this end, we first prove
a monotonicity  property of the two-dimensional embedded MC $\{{\bf X}_n\}$.
Let  ${\bf x}=(x_1,\,x_2)$ denote the point of   $\mathbb{Z}^2_+$.
Then  for all $ {\bf x}>(0,\,0)$ and $i,\,j\ge - 1$ we  define    the  probabilities
$$
\P(i,\,j):=\P\big(X^{(1)}_1=x_1+i,\,X^{(2)}_1=x_2+j|X^{(1)}_0=x_1,\,X^{(2)}_0=x_2\big),
$$
where, by definition, $\P(-1,\,-1)=0$.
Recall definition (\ref{p11}), (\ref{p22}) of the probabilities $\{p_k^{(j)}(i)\}$
and  define the total input rate to the server by $\Lambda=\lambda_1+\lambda_2+\alpha_1+\alpha_2$, when the both orbits are non-empty. We have  the following four  alternative cases implying the event
$${\cal A}=\{X^{(1)}_1=x_1+i,\, X^{(2)}_1=x_2+j\},
$$
\begin{enumerate}
    \item with probability $\lambda_1/\Lambda\,p_1^{(1)}(i)p_1^{(2)}(j),$
    when a class-1 new arrival occupies the server,
where we use the  independence of  inputs;
\item  with probability $\lambda_2/\Lambda \, p_2^{(1)}(i) p_2^{(2)}(j),$
when
a class-2 new arrival occupies the server;
\item  with probability $\alpha_1/\Lambda\,p_1^{(1)}(i+1)p_1^{(2)}(j)$,
when a class-1 orbit customer occupies the server;
\item and finally,
with probability $\alpha_1/\Lambda\, p_2^{(1)}(i)p_2^{(2)}(j+1)$, when a class-2 orbit customer occupies the server.
\end{enumerate}
Collecting together all possible cases, we obtain
\begin{eqnarray}
\P(i,\,j)&=&\frac{1}{\lambda_1+\lambda_2+\alpha_1+\alpha_2}\Big[\lambda_1p_1^{(1)}(i)p_1^{(2)}(j)+\lambda_2p_2^{(1)}(i)p_2^{(2)}(j)\nonumber \\
&+&\alpha_1p_1^{(1)}(i+1)p_1^{(2)}(j) +\alpha_2p_2^{(1)}(i)p_2^{(2)}(j+1)\Big],\label{pIJ}
\end{eqnarray}
 where   $p_1^{(k)}(-1)=0,\,p_2^{(k)}(-1)=0,\,k=1,2$, by definition. Moreover, we have
\begin{eqnarray*}
\P(i,\,-1)&=&\frac{\alpha_2}{\lambda_1+\lambda_2+\alpha_1+\alpha_2}p_2^{(1)}(i)p_2^{(2)}(0),\quad i\ge 0;\\
\P(-1,\,j)&=&\frac{\alpha_1}{\lambda_1+\lambda_2+\alpha_1+\alpha_2}p_2^{(1)}(0)p_2^{(2)}(j),\quad j\ge 0.
\end{eqnarray*}

Now, for arbitrary ${\bf y}\in \mathbb{Z}^2_+$, we define the set
$$
C_{\bf y}=\{{\bf z}\in \mathbb{Z}^2_+\,:\,{\bf z}\ge {\bf y}\}.
$$
Then, following \cite{FOSS}, we must establish  the following  monotonicity property of MC $\{{\bf X}_n\}$:
\begin{equation}\label{mono}
\P\big({\bf X}_1\in C_{\bf y} {\big |}{\bf X}_0={\bf x}\big)\ge \P\big({\bf X}_1\in C_{\bf y}\big|{\bf X}_0=\hat {\bf x}\big),\quad  {\bf x}\ge \hat {\bf x},\quad {\bf y}\in \mathbb{Z}^2_+.
\end{equation}
Next, fix arbitrary ${\bf x},\,\hat {\bf x}
>(0,\,0)$ and  calculate the left hand side of (\ref{mono}):
\begin{eqnarray}\nonumber
\P\big({\bf X}_1\in C_{\bf y}\big|{\bf X}_0={\bf x}\big)&=&\sum_{i\ge y_1-x_1}\sum_{j\ge y_2-x_2}\P\big(X_1^{(1)}=x_1+i,\,X_2^{(1)}=x_2+j\big|{\bf X}_0={\bf x}\big)\\
&=& \sum_{i\ge y_1-x_1}\sum_{j\ge y_2-x_2}\P(i,\,j).\label{left}
\end{eqnarray}
 Taking into account  $x_1\ge \hat x_1$, $x_2\ge \hat x_2$  and (\ref{pIJ}), the expression (\ref{left}) can be represented  as
\begin{eqnarray*}
   \P\big({\bf X}_1\in C_{\bf y}\big|{\bf X}_0={\bf x}\big)&=& \sum_{i\ge y_1-x_1}\big(\sum_{j= y_2-x_2}^{y_2-\hat x_2-1}\P(i,\,j)+\sum_{j\ge y_2-\hat x_2}\P(i,\,j)\big)\\
  &=& \sum_{i\ge y_1-x_1}\sum_{j= y_2-x_2}^{y_2-\hat x_2-1}\P(i,\,j)+ \sum_{i\ge y_1-x_1}\sum_{j\ge y_2-\hat x_2}\P(i,\,j)\\
  &=&
\sum_{i\ge y_1-x_1}\sum_{j= y_2-x_2}^{y_2-\hat x_2-1}\P(i,\,j)+\sum_{j\ge y_2-\hat x_2} \sum_{i= y_1-x_1}^{y_1-\hat x_1-1}\P(i,\,j) \\
&+&\P\big({\bf X}_1\in C_{\bf y}\big|{\bf X}_0=\hat {\bf x}\big),
\end{eqnarray*}
which implies the monotonicity property (\ref{mono}).

As the embedded two-dimensional MC $\{{\bf X}_n\}$ satisfies the monotonicity property and the second orbit grows to infinity by Theorem~2, we can apply
Theorem~B from the Appendix, to state the following result.

\bigskip

\noindent
{\bf Theorem 3.} {\it Let ${\bf \pi}=\{\pi_n, n\ge 0\}$ be the stationary distribution of  the orbit size in the auxiliary single-orbit system at the departure instants, and let  assumptions of Theorem~2 hold. Then,
$X_{n}^{(1)}\to X^{(1)}$
in the total variation norm where $X^{(1)}$ has distribution $\pi$.}

\bigskip

As a final remark of this section, let us explain why the probability $\P_B^{(2)}$  satisfies
expression (\ref{77}). When the auxiliary single-orbit system is stable,  $\P_B^{(2)}$  first must include  fraction  $\rho_1$  of  time when the server is occupied  by  class-1 customers. Next, the other fraction of time, $1-\rho_1$, is devoted to serving class-2 customers. When the server is working as the auxiliary  system with input rate $\lambda_2+\alpha_2$ and service rate $\mu_2$, the loss probability equals
$$
\frac{(\lambda_2+\alpha_2)/\mu_2}{1+(\lambda_2+\alpha_2)/\mu_2}=
\frac{\rho_2+\hat\rho_2}{1+\rho_2+\hat\rho_2}.
$$
Now collecting together both these fractions, we easily obtain that
$$
\rho_1+(1-\rho_1)\frac{\rho_2+\hat\rho_2}{1+\rho_2+\hat\rho_2}=\frac{\rho+\hat\rho_2}{\rho_2+\hat\rho_2+1}=\P_B^{(2)},
$$
as intuitively expected.
Note that in the analysis above we implicitly used the PASTA property allowing in our case to equate fraction of class-2 arrivals which meet server busy by other class-2 customers and the fraction of time when server is occupied by class-2 customers.

\section{Comparison with known stability results}
\label{sec:compare}

In this section, we compare  the obtained stability criterion  (\ref{Theo2}) with earlier obtained stability conditions mentioned in the introduction.  In the papers \cite{PEVA,AMNS}
the following {\it necessary} stability condition
$$
\rho=:\sum_{k=1}^K\rho_k<\min_{1\le k\le K}\frac{\alpha_k}{\lambda_k+\alpha_k},
$$
has been obtained for a bufferless $K$-class retrial system, in which class-$k$ customers follow Poisson input with rate $\lambda_k$, have i.i.d.  general
service times
with the mean $\mu_k$  and retrial rate $\alpha_k$.
As we see, for $K=2$ classes, this necessary condition coincides with stability criterion (\ref{Theo2}).   On the other hand, {\it sufficient} stability condition  from \cite{PEVA} has the form
$$
\rho<\min_k\Big(\frac{\alpha_k}{\alpha_k+\lambda}\Big),
$$
where $\lambda=\sum_{k=1}^K\lambda_k$, and definitely less tight than condition (\ref{Theo2}) (for $K=2$).

We note that for a single-class system, condition (\ref{Theo2})
becomes (in an evident notation)
\begin{eqnarray}
\rho<\frac{\alpha}{\lambda+\alpha},\label{80}
\end{eqnarray}
and coincides with  stability condition obtained in a few previous papers \cite{Fayolle1986,Lillo,MMOR,AMNS}.  Note that in \cite{Lillo} a renewal input is allowed while service time is exponential. On the other hand, system in \cite{MMOR} allows both general service time and a general renewal input.
We note that condition (\ref{80}), written as
$$
\lambda\rho<\alpha(1-\rho),
$$
has a very clear intuitive interpretation: input rate to the orbit
(generated by customers meeting busy server) must be less than the output rate from the orbit (the rate of successful attempts).
Of course a similar interpretation holds for stability conditions $\rho<\alpha_k/(\lambda_k+\alpha_k)$,
 for each orbit in the multi-class system.
One more interesting interpretation of condition (\ref{80}) is the following: when the input is Poisson,  by property PASTA,  $\rho$  is the probability that an arriving customer meets server busy and joins the orbit, while  the r.h.s of (\ref{80}) equals
\begin{eqnarray}
\frac{\alpha}{\lambda+\alpha}=\P(\xi<\tau),\nonumber
\end{eqnarray}
that is  the probability that
the retrial  time $\xi$ is less than the (remaining) interarrival time $\tau$ and thus the {\it orbital customer} occupies the server. As a result, the orbit size decreases, and this negative drift implies stability of the system.
We also note that other  related stability results
can be found in the references in the papers \cite{Questa2015,MMOR}.

Finally, we would like to mention a series of recent works devoted to regenerative stability analysis of the multiclass  retrial systems with  {\it coupled orbits} (or {\it state-dependent retrial rates}), being  a far-reaching generalization  of the constant retrial rate systems, in which the
retrial rate of each orbit depends on the binary state (busy or idle) of all other orbits, see \cite{EPEW,
Fruct23, fruct24, QTNA2019,DCCN2020, smarty}.  In particular,  this analysis is based on PASTA and a coupling procedure connecting the real processes of the retrial with the  independent Poisson processes corresponding to various `configurations' of the (binary) states of the orbits.

\section{Simulations}
\label{sec:sim}

\subsection{Exponential service times}
\label{subsec:simexp}

\begin{figure}[htb]	
\centering
\includegraphics[width=0.8\textwidth]{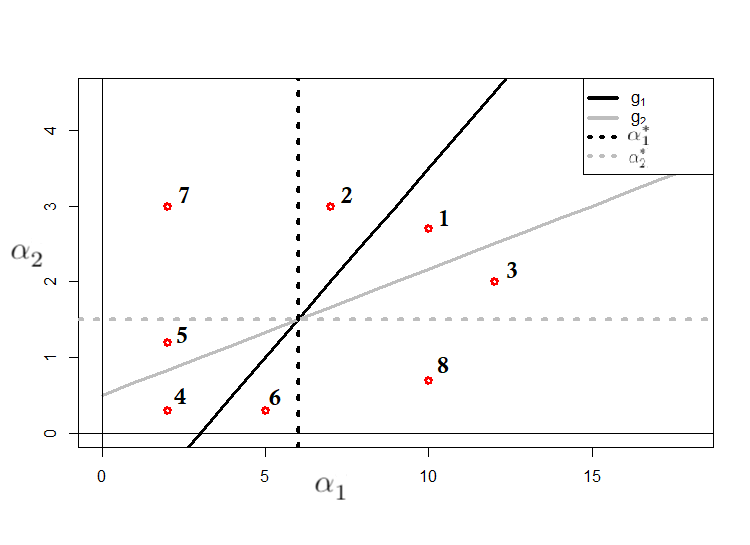}
	\caption{Stability regions for $\rho=0.75$, exponential service times}
\label{sim1}
\end{figure}

First we consider the case of exponential service times with corresponding service rates $\mu_1,\,\mu_2$ and consider a particular case $\lambda_1=2,\,\lambda_2=0.5,\,\mu_1=4,\,\mu_2=2$. Thus,
$$
\rho_1=0.5,\,\rho_2=0.25,\,\rho=0.75.
$$
Stability regions for such a configuration are presented in Figure~\ref{sim1}.
The regions in Figure~\ref{sim1} correspond to various cases of Theorems 1 and A. Namely,
\begin{itemize}
    \item[1.] (a.1)--stability;
    \item[2]  (b.1)--stability;
    \item[3.] (c.1)--stability;
    \item[4.] (d.1)--transience. 
\end{itemize}
Note that  (b.2)--transience region is defined by conditions
\begin{eqnarray*}
\alpha_2&\ge&g_1,\\
\alpha_2&>&g_2,\\
\alpha_1&<& \alpha_1^*.
\end{eqnarray*}
Figure~\ref{sim1} illustrates that such a region  
is divided into two subsets: for $\alpha_2<\alpha_2^*$ 
and  $\alpha_2>\alpha_2^*$. In terms of load coefficients, the condition $\alpha_2<\alpha_2^*$ transforms to $\rho>  \alpha_2/(\lambda_2+\alpha_2)$. Thus, we define two separate cases as follows:
\begin{itemize}
    \item[5.]  (b.2)--transience and $\rho>  \alpha_2/(\lambda_2+\alpha_2)$;
    \item[7.]  (b.2)--transience and  $\rho<  \alpha_2/(\lambda_2+\alpha_2)$.
    \end{itemize}
The same phenomenon arises in    (c.2)--transience  region, which is divided by the horizontal line $\alpha_2=\alpha_2^*$. Then we also include in the simulation plan the cases
\begin{itemize}
    \item[6.]  (c.2)--transience and  $\rho>  \alpha_1/(\lambda_1+\alpha_1)$;
    \item[8.]  (c.2)--transience and  $\rho<  \alpha_1/(\lambda_1+\alpha_1)$.
    \end{itemize}

We recall  that by
Theorem 2 the   cases 5 and  7 (and  also 6 and  8) correspond to locally  stable scenario: one orbit is tight while the other orbit increases to infinity in probability.
%


Simulation results for particular cases, which correspond to all possible stability/instability regions, are presented in Table~\ref{tab1}.
Note that for considered configurations we obtain
$$\alpha^*_1=6,\,\alpha^*_2=1.5,$$
and we recall that the conditions $\rho<\alpha_k/(\lambda_k+\alpha_k)$ are equivalent to the bounds $\alpha_k>\alpha^*_k$, $k=1,2$.
The last $2$ columns in Table~\ref{tab1} are based on simulation results. The mark ``yes'' means that the orbit has stable dynamics, while the mark ``no'' indicates the growth to infinity.

\begin{table}
\centering
\caption{Simulation results for the system with exponential service times,  $\rho=0.75$.}
\label{tab1}
\begin{tabular}{|l|l|l|l|l|l|l|l|l|}
\hline
№ &  $\alpha_1$ & $\alpha_2$&$\frac{\alpha_1}{\lambda_1+\alpha_1}$&$\frac{\alpha_2}{\lambda_2+\alpha_2}$&$g_1(\alpha_1)$&$g_2(\alpha_1)$&stable  $X^{(1)}$&stable $X^{(2)}$\\
\hline
1.&10.00 &2.70 &0.83&0.84&3.50 &2.17 &yes &yes\\
2.& 7.00&3.00 &0.78&0.86& 2.00& 1.17& yes&yes\\
3.&12.00 &2.00 &0.86&0.80& 4.50& 2.50& yes&yes\\
4.&2.00 &0.30 &0.50&0.38&-0.50&0.83 &no &no\\
5.&2.00 &1.20 &0.50&0.71&-1.00 &1.67 & no & yes\\
6.& 5.00&0.30 &0.71&0.38& 1.00 &1.33 &yes&no\\
7.&2.00 &3.00 &0.50 &0.86&-0.50 & 0.83& no&yes\\
8.& 10.00& 0.70&0.83&0.58& 3.50& 2.17& yes& no\\
\hline
\end{tabular}
\end{table}

Region 7 
corresponds to the set, where the condition  $\rho<\alpha_2/(\alpha_2+\lambda_2)$ holds and the condition $\rho<\alpha_1/(\alpha_1+\lambda_1)$ is violated. In this region, as expected from the theoretical results, we obtained that only the second orbit is stable. The symmetric result was obtained for region 8: we have  $\rho<\alpha_1/(\alpha_1+\lambda_1)$, $\rho>\alpha_2/(\alpha_2+\lambda_2)$ and consequently only the first orbit is stable. Thus, at the first sight, it appears that the condition $\rho<\alpha_k/(\alpha_k+\lambda_k)$ defines the (local) stability of $k$-class orbit.

It is rather surprising what we observe as simulations results in cases 5 and 6.
In these regions both conditions 
are violated, while the $2$-nd (the $1$-st)  orbit is stable  in case 5 (6). Note that regions 5 and 7 correspond to the transient case (b.2) from Theorem~A, while regions 6 and 8 correspond to the transient case (c.2). The only case when the both orbits are unstable, was obtained in region 4, which corresponds to transient case (d) from Theorem~A and is defined by $g_1(\alpha_1)<\alpha_2<g_2(\alpha_1)$.

\begin{figure}[htb]	
\centering
\includegraphics[width=0.7\textwidth]{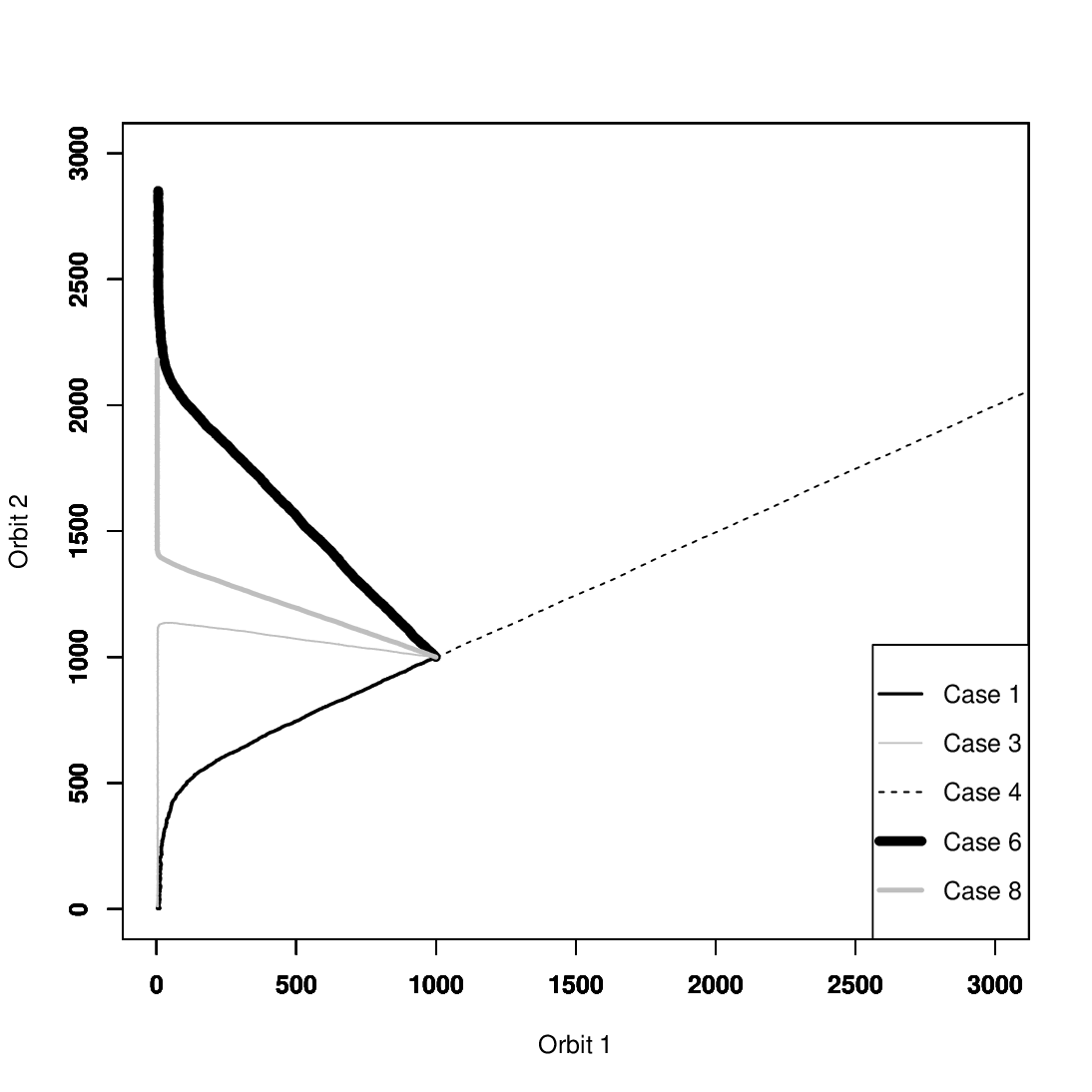}
	\caption{Averaged orbits' dynamics for non-zero initial state, exponential service times
}
\label{figjm1}
\end{figure}

Next, we consider non-zero initial states of both orbits:
$$X^{(1)}(t_1^-)=X^{(2)}(t_1^-)=1000,$$
and explore the orbit behavior, setting values $(\alpha_1,\,\alpha_2)$ from Table~\ref{tab1}. The simulation results are presented as a 2D plot in Figure~\ref{figjm1}.  Note that cases 2, 5 and 7 are symmetric to 3, 6 and 8, respectively. Orbits' dynamics (stable/unstable) from Figure~\ref{figjm1} correspond to the results for those of the zero-initial state system (see Table \ref{tab1}). Both cases 6 and 8 illustrate the phenomenon of partial stability (only the second orbit grows, while the first unloads). Note that the bold black line, corresponding to case 6, majorizes the bold grey line, which describes the configuration 8. Such results are explained by the fact that unlike case 8, in case 6 the condition $\rho<\alpha_1/(\lambda_1+\alpha_1)$ is violated.

\begin{figure}[htb]	
\centering
\includegraphics[width=0.7\textwidth]{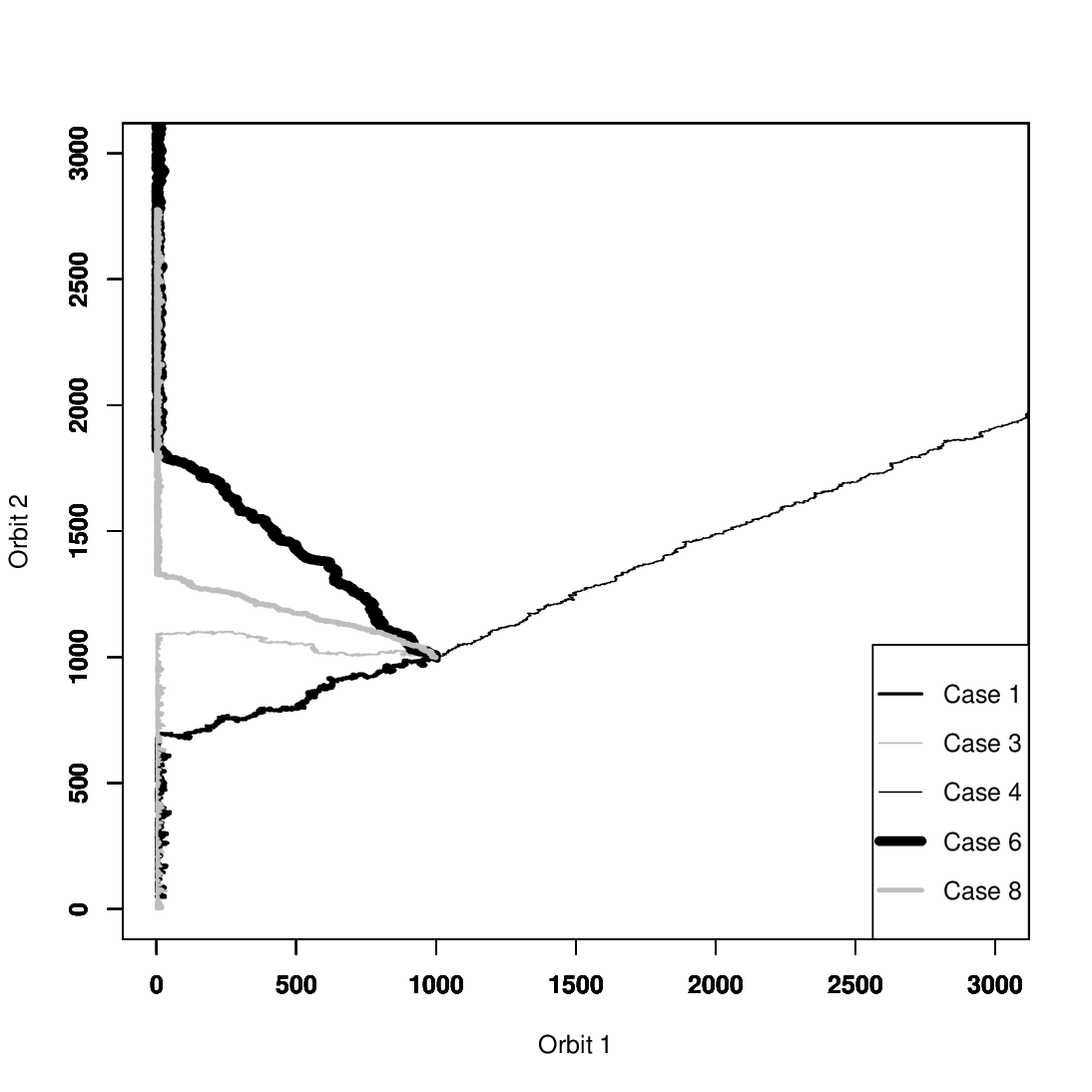}
	\caption{ Dynamics of the orbits based on  a single  trajectory: non-zero initial state and  exponential service times.}
\label{figjm1-1}
\end{figure}

Figure~\ref{figjm1} depicts the empirical means obtained by averaging over $m=100$ independent trajectories,  while Figure~\ref{figjm1-1} 
demonstrates the corresponding results based on only one realisation, and it is seen
that the orbit dynamics in both figures are in agreement.

\subsection{Pareto service times}
\label{subsec:simPareto}

Now we consider Pareto class-dependent distributions of service times with shape parameter $x_k>0$ and degree parameter $\beta_k>0,\,k=1,\,2$. Namely,
$$
\P(S^{(k)}\le x)=1-(x_k/x)^{\beta_k},\qquad x>x_k
$$
and, consequently,
$$
\E S^{(k)}=\frac{x_k\beta_k}{\beta_k-1}.
$$
Then we set again $\lambda_1=2,\,\lambda_2=0.5$ and set $x_1=0.125,$ $x_2=0.4$, $\beta_1=2$, $\beta_2=5$, which yields
$\E S^{(1)}=0.25,\, \E S^{(2)}=0.5$. Thus, similarly to the case of exponential service times,  we obtain
$$
\rho_1=0.5,\quad \rho_2=0.25,\quad \rho=0.75.
$$
Next we define the  values of retrial rates $(\alpha_1,\,\alpha_2)$, corresponding to all cases presented in Table~\ref{tab1}.

\begin{figure}[htb]	
\centering
\includegraphics[width=0.7\textwidth]{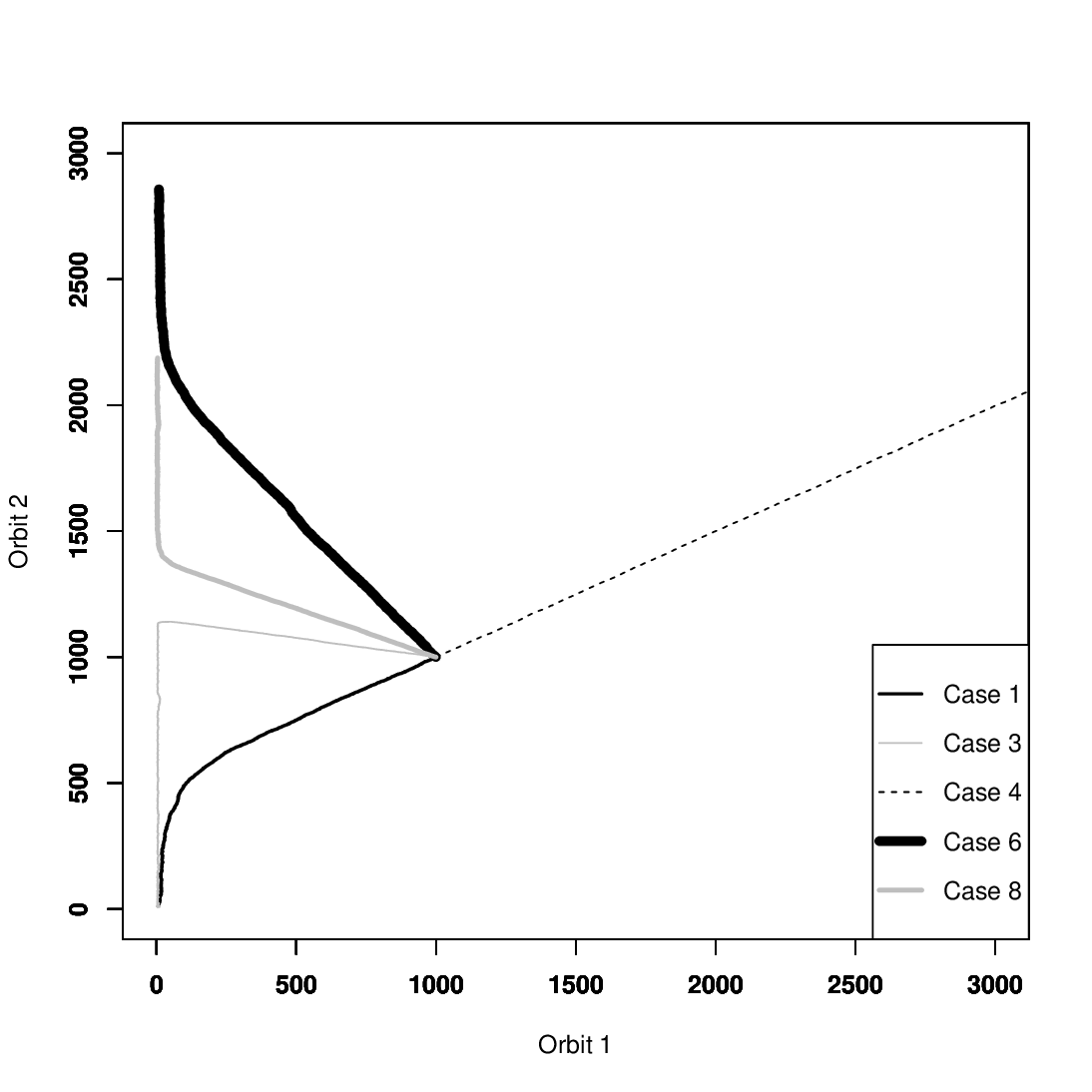}
	\caption{Averaged orbits' dynamics for non-zero initial state, Pareto service times}
\label{figgjpar}
\end{figure}

Simulations for the model with Pareto distribution of service times and  non-zero initial conditions: $X^{(1)}(t_1^-)=X^{(2)}(t_1^-)=1000$ are illustrated in Figure~\ref{figgjpar} and follow closely the results for exponential service times.
It is worth mentioning that
simulation results illustrate the  same  phenomenon of partial stability in cases 6 and 8 as we detected for exponential service times although  Pareto distribution {\it does not belong to class $\mathcal{D}$}.
It shows that the latter requirement is rather technical and the statement of Theorem~2 should hold for
a wider class of distribution functions.

\newpage
\section*{Appendix A}
In this appendix we present known results which are used in the main part of the paper.

First we consider two-dimensional MC $\{{\bf X}_n\}$ and mention the conditions for applicability of a basic theorem from \cite{FMM}.

\bigskip

\noindent
{\bf Condition A. ({\it lower boundedness condition})} \\
 \begin{equation*}
\left\{
 \begin{gathered}
\P\Big(X^{(1)}_{n+1}=i,X^{(2)}_{n+1}=j \Big|X^{(1)}_n>0,\,X^{(2)}_n>0\Big)=0,\qquad \text{if }i<-1,\qquad \text{or }j<-1, \hfill
 \\
\P\Big(X^{(1)}_{n+1}=i,X^{(2)}_{n+1}=j \Big|X^{(1)}_n>0,\,X^{(2)}_n=0\Big)=0,\qquad \text{if }i<-1,\qquad \text{or }j<0,\hfill
 \\
\P\Big(X^{(1)}_{n+1}=i,X^{(2)}_{n+1}=j \Big|X^{(1)}_n=0,\,X^{(2)}_n>0\Big)=0,\qquad \text{if }i<0,\qquad \text{or }j<-1. \hfill
 \\
 \end{gathered}
 \right. \label{condA}
 \end{equation*}
Note that the definitions of the transition probabilities in Section~2.1 automatically imply
the validity of Condition~A for our embedded MC.

Recall that $\Delta_{n+1}={\bf X}_{n+1}-{\bf X}_n$.

\bigskip

\noindent
{\bf Condition B. ({\it first moment condition})}
  \begin{equation}
  \E\Big[||\Delta_{n+1}||\Big|{\bf X}_n=(k,\,l)\Big]\le C<\infty,\qquad \forall (k,\,l)\in \mathbb{Z}_+^2,\label{condB}
  \end{equation}
where $||\cdot||$  denotes the Euclidean norm and $C>0$ is a constant.

Since
\begin{eqnarray*}
\E\Big[ ||\Delta_{n+1}||\Big]&\le& \E[\Delta_{n+1}^{(1)}+\Delta_{n+1}^{(2)}]\\
&\le& 2\max(\lambda_1,\,\lambda_2,\,\alpha_1,\,\alpha_2)\max(\E S^{(1)},\,\E S^{(2)})<\infty,
\end{eqnarray*}
the condition (\ref{condB}) holds for all sets of (finite) parameters for our MC.

Now recall that   $\M^{01}_k,\,\M^{10}_k,\,\M^{11}_k$ denote the mean increments (drifts)
of the $k$-th component of the MC $\{{\bf  X}_n\}$ between the  $n$-th and  $(n+1)$-st departure instants, given the respective conditions $(X_n^{(1)}=0,\,X_n^{(2)}>0),\,(X_n^{(1)}>0,\,X_n^{(2)}=0)$ or $(X_n^{(1)}>0,\,X_n^{(2)}>0)$.

\bigskip

\noindent
{\bf Theorem A. (See Theorem~3.3.1 in \cite{FMM})} Let a MC $\{{\bf  X}_n\}$ be aperiodic and irreducible and let Conditions A and B hold.\\
{\bf (a)} If $\M_1^{11}<0,\,\M_2^{11}<0$, then MC $\{{\bf  X}_n\}$ is
\begin{enumerate}
    \item ergodic (positive recurrent) if
    \begin{equation}\label{first}
\left\{
\begin{gathered}
\M_1^{11}\M_2^{10}-\M_2^{11}\M_1^{10}<0, \hfill
\\
\M_2^{11}\M_1^{01}-\M_1^{11}\M_2^{01}<0,\hfill
\\
\end{gathered}
\right.
\end{equation}
 \item non-ergodic 
 if either
    \begin{equation}\label{nonerg}
        \M_1^{11}\M_2^{10}-\M_2^{11}\M_1^{10}\ge 0\,\,\,\,\text{ or }\,\,\,\, \M_2^{11}\M_1^{01}-\M_1^{11}\M_2^{01}\ge0.
    \end{equation}
\end{enumerate} 
{\bf (b)} If $\M_1^{11}\ge 0,\,\M_2^{11}<0$, then MC $\{{\bf  X}_n\}$ is
\begin{enumerate}
    \item  ergodic if
 \begin{equation}\label{first_b}
\M_1^{11}\M_2^{10}-\M_2^{11}\M_1^{10}< 0,
\end{equation}
 \item transient if
  \begin{equation} \label{trb}
\M_1^{11}\M_2^{10}-\M_2^{11}\M_1^{10}>0.
\end{equation}
\end{enumerate}
{\bf (c)} If $\M_1^{11}< 0,\,\M_2^{11}\ge0$, then MC $\{{\bf  X}_n\}$ is
\begin{enumerate}
    \item  ergodic if
 \begin{equation}\label{first_c}
\M_2^{11}\M_1^{01}-\M_1^{11}\M_2^{01}< 0,
\end{equation}
 \item transient if
  \begin{equation}\label{trc}
\M_2^{11}\M_1^{01}-\M_1^{11}\M_2^{01}>0.
\end{equation}
\end{enumerate}
{\bf (d)} If $\M_1^{11}\ge 0,\,\M_2^{11}\ge0,\,\M_1^{11}+\M_2^{11}>0$, then MC $\{{\bf X}_n\}$ is transient.

\medskip


Next we present partial stability results from \cite{FOSS} for a two-dimensional MC.

\bigskip

\noindent
{\bf Theorem B. (See Proposition~2 in \cite{FOSS})} {\it For a MC ${\bf X}_n$, with a state space $\mathbb{Z}^2_+$, assume the following:
\begin{enumerate}
    \item $X_n^{(2)}\to \infty$ in probability, as $n\to \infty$, given $X_0^{(2)}=0,\,X_0^{(1)}=0$;
    \item $\P\big( X_{n+1}^{(1)}=j\big| X_{n}^{(1)}=i,\, X_{n}^{(2)}=l\big)=p_{ij}$ for all values of $l>0$, where
$p_{ij}$ are transition probabilities of an ergodic MC with the unique stationary distribution $\pi=\{\pi_j\}$;
\item MC  ${\bf X}_n$ is monotone. Namely,
$$
\P\big({\bf X}_1\in C_{\bf y} {\big |}{\bf X}_0={\bf x}\big)\ge \P\big({\bf X}_1\in C_{\bf y}\big|{\bf X}_0=\hat {\bf x}\big),\quad  {\bf x}\ge \hat {\bf x},\quad {\bf y}\in \mathbb{Z}^2_+.
$$
\end{enumerate}
Then, for all initial states $i_0,\,j_0$:
\begin{equation}
     \sup_j \big|\P\big(X_{n}^{(1)}=j| X_{0}^{(1)}=i_0,\, X_{0}^{(2)}=j_0\big)-\pi_j\big|\to 0, \quad n\to \infty,
\end{equation}
i.e., $X_{n}^{(1)}$ converges in distribution to $\pi$ in the total variation norm.}

\section*{Appendix B}

Below we present some technical details of the proof of Theorem 1.
Note that our goal is to apply the results of Theorem A to the system under consideration. First, consider case (a) from Theorem A and recall  expressions (\ref{beMM1})--(\ref{beM2}) for $\M_1^{10},\,\M_2^{10},\,\M_1^{11},\,\M_2^{11}$, respectively.
Conditions $\M^{11}_1<0,\,\M^{11}_2<0$ can be relatively easy rewritten as (\ref{IdM1}) and (\ref{IdM2}) in terms of the load coefficients. Next let us also represent the system (\ref{first})
in terms of the load coefficients. Define the following auxiliary
parameters
\begin{eqnarray*}
A_1&=&(\lambda_1+\alpha_1)\lambda_1\/\mu_1+\lambda_2\lambda_1/\mu_2,\\
A_2&=&(\lambda_1+\alpha_1)\lambda_2/\mu_1+\lambda_2\lambda_2/\mu_2,\\
A_3&=&\lambda_1+\lambda_2+\alpha_1,\\
A_4&=&\lambda_1+\lambda_2+\alpha_1+\alpha_2.
\end{eqnarray*}
Thus, from (\ref{beMM1})--(\ref{beM2}) we obtain
\begin{eqnarray*}
\M^{10}_1&=&\frac{A_1} {A_3}-\frac{\alpha_1}{A_3},\\
\M^{10}_2&=&\frac{A_2} {A_3},\\
\M_1^{11}&=&\frac{A_1}{A_4}+\frac{\alpha_2\lambda_1/\mu_2-\alpha_1}{A_4},\\
\M_2^{11}&=&\frac{A_2}{A_4}+\frac{\alpha_2\lambda_2/\mu_2-\alpha_2}{A_4}.
\end{eqnarray*}
Next, denote
\begin{eqnarray*}
B_1&=&\frac{\alpha_2\lambda_1/\mu_2-\alpha_1}{A_4},\\
B_2&=&\frac{\alpha_2\lambda_2/\mu_2-\alpha_2}{A_4},\\
B_3&=&\frac{\alpha_1}{A_3}.
\end{eqnarray*}
Thus, the condition $\M_1^{11}\M_2^{10}-\M_2^{11}\M_1^{10}<0$ 
transforms to
\begin{equation*}
    \Big(\frac{A_1}{A_4}+B_1\Big)\frac{A_2}{A_3} -\Big(\frac{A_2}{A_4}+B_2\Big)\Big(\frac{A_1} {A_3}-B_3\Big)<0.
\end{equation*}
After opening brackets, we obtain
\begin{equation*}
   B_1\frac{A_2} {A_3} +B_3\frac{A_2} {A_4}-B_2\frac{A_1}{A_3}+B_2B_3<0.
\end{equation*}
Then we substitute back the expressions for $B_1,\,B_2,\,B_3$ to get
\begin{equation}\label{calc}
   \frac{(\alpha_2\lambda_1/\mu_2-\alpha_1)A_2}{ A_3A_4} +\frac{\alpha_1A_2}{A_3A_4}-\frac{(\alpha_2\lambda_2/\mu_2-\alpha_2)A_1}{A_3A_4}+\frac{(\alpha_2\lambda_2/\mu_2-\alpha_2)\alpha_1}{A_3A_4}<0.
\end{equation}
Multiplying both sides of (\ref{calc}) by $A_3A_4>0$, yields
\begin{equation*}
    A_2\alpha_2\lambda_1/\mu_2-A_1\alpha_2(\lambda_2/\mu_2-1)+\alpha_1\alpha_2(\lambda_2/\mu_2-1)<0.
\end{equation*}
Then, after substituting the expressions for $A_1$ and $A_2$,
we obtain the following inequality in terms of the original parameters
\begin{equation*}
(\lambda_1+\alpha_1)\lambda_1/\mu_1+\lambda_2\lambda_1/\mu_2+\alpha_1(\lambda_2/\mu_2-1)<0,
\end{equation*}
which easily transforms to
\begin{eqnarray*}
(\lambda_1+\alpha_1)(\lambda_1/\mu_1+\lambda_2/\mu_2)<\alpha_1.
\end{eqnarray*}
Thus, we obtain the condition
\begin{equation*}
\rho<\frac{\alpha_1}{\lambda_1+\alpha_1}\equiv \frac{\hat \rho_1}{\rho_1+\hat\rho_1}.
\end{equation*}
After similar derivations, the condition $\M_2^{11}\M_1^{01}-\M_1^{11}\M_2^{01}<0$ is
transformed to
\begin{equation*}
\rho<\frac{\hat \rho_2}{\rho_2+\hat\rho_2}
.
\end{equation*}
Thus, the system (\ref{first}) is equivalent to $\rho<\hat\rho_k/(\hat \rho_k+\rho_k),\, k=1,2$.

\bibliographystyle{plain}
\bibliography{mybibfile}

\begin{thebibliography}{10}

\bibitem{FOSS}
I.~Adan, S.~Foss, S.~Shneer, and G.~Weiss.
\newblock Local stability in a transient {M}arkov chain.
\newblock {\em Statistics and Probability Letters}, 165(108855):1--6, 2020.

\bibitem{Artalejo}
J.R. Artalejo.
\newblock Accessible bibliography on retrial queues.
\newblock {\em Mathematical and computer modelling}, 30(3-4):1--6, 1999.

\bibitem{Artalejo1}
J.R. Artalejo.
\newblock Accessible bibliography on retrial queues: progress in 2000--2009.
\newblock {\em Mathematical and computer modelling}, 51(9-10):1071--1081, 2010.

\bibitem{Artalejo2008book}
J.R. Artalejo and A.~G{\'o}mez-Corral.
\newblock {\em Retrial Queueing Systems: {A} Computational Approach}.
\newblock Springer, 2008.

\bibitem{Asmus}
S.~Asmussen.
\newblock {\em Applied Probability and Queues. 2nd edn.}
\newblock Springer, New York, 2003.

\bibitem{MMOR}
K.~Avrachenkov and E.~Morozov.
\newblock Stability analysis of {GI}/{GI}/c/{K} retrial queue with constant
  retrial rate.
\newblock {\em Mathematical Methods of Operations Research}, 79(3):273--291,
  2014.

\bibitem{AMNS}
K.~Avrachenkov, E.~Morozov, R.~Nekrasova, and B.~Steyaert.
\newblock Stability analysis and simulation of {N}-class retrial system with
  constant retrial rates and {P}oisson inputs.
\newblock {\em Asia-Pacific Journal of Operational Research}, 31(02):1440002,
  2014.

\bibitem{Questa2015}
K.~Avrachenkov, E.~Morozov, and B.~Steyaert.
\newblock Sufficient stability conditions for multi-class constant retrial rate
  systems.
\newblock {\em Queueing Systems}, 82(1-2):149--171, 2016.

\bibitem{Nain}
K.~Avrachenkov, P.~Nain, and U.~Yechiali.
\newblock A retrial system with two input streams and two orbit queues.
\newblock {\em Queueing Systems}, 77(1):1--31, 2014.

\bibitem{PEISarticle}
K.~Avrachenkov and U.~Yechiali.
\newblock Retrial networks with finite buffers and their application to
  internet data traffic.
\newblock {\em Probability in the Engineering and Informational Sciences},
  22(4):519--536, 2008.

\bibitem{KostiaUri}
K.~Avrachenkov and U.~Yechiali.
\newblock On tandem blocking queues with a common retrial queue.
\newblock {\em Computers \& Operations Research}, 37(7):1174--1180, 2010.

\bibitem{dim3}
I.~Dimitriou.
\newblock Modeling and analysis of a relay-assisted cooperative cognitive
  network.
\newblock {\em Proceedings Analytical and Stochastic Modelling Techniques and
  Applications (ASMTA)}, pages 47--62, 2017.

\bibitem{dimi}
I.~Dimitriou.
\newblock A queueing system for modeling cooperative wireless networks with
  coupled relay nodes and synchronized packet arrivals.
\newblock {\em Performance Evaluation}, 114:16--31, 2017.

\bibitem{dim2}
I.~Dimitriou.
\newblock A two-class retrial system with coupled orbit queues.
\newblock {\em Probability in the Engineering and Informational Sciences},
  31(2):139--179, 2017.

\bibitem{Dimitriou2018EJOR}
I.~Dimitriou.
\newblock A two-class queueing system with constant retrial policy and general
  class dependent service times.
\newblock {\em European Journal of Operational Research}, 270(3):1063--1073,
  2018.

\bibitem{dimitriou2019power}
I.~Dimitriou.
\newblock On the power series approximations of a structured batch arrival
  two-class retrial system with weighted fair orbit queues.
\newblock {\em Performance Evaluation}, 132:38--56, 2019.

\bibitem{Dimitriou2020Valuetools}
I.~Dimitriou and T.~Phung-Duc.
\newblock Analysis of cognitive radio networks with cooperative communication.
\newblock In {\em Proceedings of the 13th EAI International Conference on
  Performance Evaluation Methodologies and Tools (ValueTools)}, pages 192--195,
  2020.

\bibitem{FalinSurvey}
G.~Falin.
\newblock A survey of retrial queues.
\newblock {\em Queueing systems}, 7(2):127--167, 1990.

\bibitem{Falin1997}
G.~Falin and J.G.C. Templeton.
\newblock {\em Retrial queues}, volume~75.
\newblock CRC Press, 1997.

\bibitem{Fayolle1986}
G.~Fayolle.
\newblock A simple telephone exchange with delayed feedbacks.
\newblock In {\em Proc. of the international seminar on Teletraffic analysis
  and computer performance evaluation}, pages 245--253, 1986.

\bibitem{FMM}
G.~Fayolle, V.A. Malyshev, and M.V. Menshikov.
\newblock {\em Topics in the Constructive Theory of Countable Markov Chains.
  1st edn.}
\newblock Cambridge University Press, 1995.

\bibitem{Lillo}
R.E. Lillo.
\newblock A {G}/{M}/1-queue with exponential retrial.
\newblock {\em {TOP}}, 4:99--–120, 1996.

\bibitem{MorozovDelgado}
E.~Morozov and R.~Delgado.
\newblock Stability analysis of regenerative queues.
\newblock {\em Automation and remote control}, pages 1977--1991, 2009.

\bibitem{EPEW}
E.~Morozov and I.~Dimitriou.
\newblock Stability analysis of a multiclass retrial system with coupled orbit
  queues.
\newblock {\em Proceedings of 14th European Workshop, EPEW}, 10497:85–98,
  2017.

\bibitem{Fruct23}
E.~Morozov and T.~Morozova.
\newblock Analysis of a generalized retrial system with coupled orbits.
\newblock {\em Proceeding 23rd Conference of Open Innovations Association
  (FRUCT)}, pages 253--260, 2018.

\bibitem{QTNA2019}
E.~Morozov and T.~Morozova.
\newblock A coupling-based analysis of a multiclass retrial system with
  state-dependent retrial rates.
\newblock {\em Proceedings 14th International Conference on Queueing Theory and
  Network Applications}, 11688:34--50, 2019.

\bibitem{DCCN2020}
E.~Morozov and T.~Morozova.
\newblock The remaining busy time in a retrial system with unreliable servers.
\newblock {\em Proceedings International Conference on Distributed Computer and
  Communication Networks}, 12563:555--566, 2020.

\bibitem{smarty}
E.~Morozov, T.~Morozova, and I.~Dimitriou.
\newblock Simulation of multiclass retrial system with coupled orbits.
\newblock {\em Proceedings of the First International Workshop on Stochastic
  Modeling and Applied Research of Technology Petrozavodsk}, pages 6--16, 2018.

\bibitem{fruct24}
E.~Morozov, T.~Morozova, and I.~Dimitriou.
\newblock A multiclass retrial system with coupled orbits and service
  interruptions: Verification of stability conditions.
\newblock {\em Proceedings of the 24th Conference of Open Innovations
  Association}, 24:75--81, 2019.

\bibitem{Tuan}
E.~Morozov and T.~Phung-Duc.
\newblock Regenerative analysis of two-way communication orbit-queue with
  general service time.
\newblock {\em Proceedings International Conference Queueing Theory and Network
  Applications}, 10932:22--32, 2018.

\bibitem{PEVA}
E.~Morozov, A.~Rumyantsev, S.~Dey, and T.G. Deepak.
\newblock Performance analysis and stability of multiclass orbit queue with
  constant retrial rates and balking.
\newblock {\em Performance Evaluation}, 134:102005, 2019.

\bibitem{MSbook2021}
E.~Morozov and B.~Steyaert.
\newblock {\em Stability Analysis of Regenerative Queueing Models: Mathematical
  Methods and Applications}.
\newblock Springer, 2021.

\bibitem{PhungDuc2016}
T.~Phung-Duc, W.~Rogiest, Y.~Takahashi, and H.~Bruneel.
\newblock Retrial queues with balanced call blending: analysis of single-server
  and multiserver model.
\newblock {\em Annals of Operations Research}, 239(2):429--449, 2016.

\end{thebibliography}
\end{spacing}
\end{document}